# Groups with infinite linearly ordered products

by Vincent Bagayoko
IMJ-PRG (Paris)

**Abstract**

We introduce a formalism of infinite, linearly ordered products in general groups. Using this, we define infinite compositions in certain groups of formal power series such as transseries. We show that such groups can sometimes be represented as infinite, linearly ordered, semidirect products of ordered Abelian groups.

## Introduction

A linearly ordered set $I = (\underline{I}, <)$ induces an $I$-indexed product function $\Pi_I^{\min}$ on a group $(\mathcal{G}, \cdot, 1)$, whose domain is the set $\mathcal{G}^{(\underline{I})}$ of finitely supported functions $\underline{I} \longrightarrow \mathcal{G}$, and which sends a $g \in \mathcal{G}^{(\underline{I})}$ to the finite product of its values in order, i.e.

$$\Pi_I^{\min} g = g(i_1) \cdots g(i_n), \tag{1}$$

where the support of $g$ is $\{i_1, \ldots, i_n\}$ and $i_1 < \cdots < i_n$. How to extend (1) to families $\underline{I} \longrightarrow \mathcal{G}$ that are not finitely supported?

We propose sensible generalisations of some of the properties of the family $(\Pi_I^{\min})_I$ on a general group $\mathcal{G}$, where $I$ ranges among all linearly ordered sets. Considering a family of partial functions $\Pi_I : \mathcal{G}^{\underline{I}} \longrightarrow \mathcal{G}$ indexed by linearly ordered sets $I = (\underline{I}, <)$, we give axioms pertaining to the existence and values of possibly infinite products

$$\Pi_I g = \cdots g(i) \cdots g(j) \cdots \in \mathcal{G},$$

where $i, j \in I$ are thought as occurring according to the ordering $i < j$, and $g : \underline{I} \longrightarrow \mathcal{G}$ is a function. The group $(\mathcal{G}, \cdot, 1)$ together with the family $(\Pi_I)_I$ indexed by *all* linearly ordered sets $I$, is called a multipliability group.

Our motivation for studying infinite products comes from the theory of linearly ordered groups. We introduced [3] an elementary class of ordered groups, called “growth order groups”, that includes certain linearly ordered groups of generalised power series under formal composition laws [9]. We believe such groups of formal series have natural structures of multipliability groups. We also expect that the formal Lie correspondence of [4] should extend in most generality to sufficiently rich multipliability groups on the Lie group side, and summability algebras with evaluations on the Lie algebra side (see Remark 4.17). In this paper, we illustrate these aspects for the field $\mathbb{S}$ of *flat transseries* (see [1, Section 3]), a subfield of Dahn-Göring [5] and Ecalle's [7] field of log-exp transseries.

Dating to Hahn [9], Hausdorff [10], Krull [15] and Kaplansky [14], the theories of extensions of (ordered) algebraic structures benefited from both the development of valuation theory and the construction of goups and fields of formal series. In order to study the non-commutative, compositional context, we require an extension of the theory of non-commutative valuations on growth order groups (see [3, Section 2.2]), a general way to construct (ordered) groups of formal non-commutative series, and a way to study the properties of infinite non-commutative expansions and infinite products in general. This paper is dedicated to the third task, with some applications to the second one.

Looking at infinite products in the context of groups of formal series under composition laws, we are to develop notions of transfinite ordered compositions. In the case of certain types of logarithmic transseries, transfinite composition appear when studying normalisation problems related to the dynamics of (germs of) diffeomorphisms of the complex variable which fix the origin [19, 1], as well as compositional inversion problems [2]. Here, we provide a general theory for such transfinite compositions, and we will give a thorough study of the case of flat transseries. We will rely on it in future work in order to study and solve equations and inequalities over certain ordered groups.

Let us describe our main results. Section 1 introduces our conventions for groups of functions and linearly ordered sets. In Section 2, we give the axiomatic definition of multipliability groups (axioms **MG1**–**MG7**), adapting the definition in [4, Section 1] to the non-commutative context.

In order to construct the most important examples of multipliability groups, we rely on abstract axiomatic notions of infinite sums in vector spaces as considered in [23, 8, 4]. In [4], associative algebras $A$ over a field $k$ are enriched with an infinitary operator $\Sigma$ of transfinite summation, where each set $I$ gives rise to a linear function $\Sigma_I : \operatorname{dom} \Sigma_I \longrightarrow A$ whose domain is a subspace of $A^I$, and which takes a so-called summable family $a \in \operatorname{dom} \Sigma_I$ to its sum $\Sigma_I\, a \in A$. An example is the algebra $k\langle\langle I \rangle\rangle$ of formal power series in a set $\{X_i : i \in I\}$ of non-commuting variables. If moreover $A = k + \mathfrak{m}$ is a local algebra with maximal ideal $\mathfrak{m}$, then we say that $A$ *has evaluations* if for all summable families $a : I \longrightarrow \mathfrak{m}$, the family $(a(i_1) \cdots a(i_n))_{n \in \mathbb{N} \wedge (i_1, \ldots, i_n) \in I^n}$ is also summable. Then each summable $a$ gives induces an evaluation morphism $\operatorname{ev}_a : k\langle\langle I \rangle\rangle \longrightarrow A$ which sends $X_i$ to $a(i)$ for each $i \in I$, and preserves infinite sums. All these notions are presented in Section 3.

In Section 4, we define ordered products in groups of the form $\mathcal{G} = 1 + \mathfrak{m}$ where the algebra $A = k + \mathfrak{m}$ has evaluations. For a linearly ordered set $(I, <)$ the natural ordered product $\cdots (1 + X_{i'}) \cdots (1 + X_i) \cdots (1 + X_{i''}) \cdots$ where the terms appear in order $i' < i < i''$, is defined to be the series

$$\mathrm{OP}_I := \sum_{n \in \mathbb{N} \wedge i_1 < \cdots < i_n} X_{i_1} \cdots X_{i_n} \in k\langle\langle I \rangle\rangle$$

This yields a notion of ordered products in the subgroup $1 + \mathfrak{m} \subseteq A^{\times}$: if $a : I \longrightarrow \mathfrak{m}$ is summable, then $1 + a : I \longrightarrow 1 + \mathfrak{m}$ is said multipliable, and the ordered product $\Pi_I(1 + a)$ is defined as the evaluation $\operatorname{ev}_a(\mathrm{OP}_I)$ of $\mathrm{OP}_I$ at $a$.

**Theorem 1.** [Theorem 4.15] *Let $(A, \Sigma)$ be a summability algebra with evaluations, with maximal ideal $\mathfrak{m}$. For each linearly ordered set $I$, let $\Pi_I$ be the ordered product function on $(1 + \mathfrak{m})^{\underline{I}}$. Then $(1 + \mathfrak{m}, 1, \cdot, \Pi)$ is a multipliability group.*

The group $1\text{-}\mathrm{Aut}_k^+(\mathbb{K})$ of linear automorphisms $\sigma$ of a valued field $(\mathbb{K}, v)$ of Hahn series which commute with infinite sums and satisfy $v(\sigma(a) - a) > v(a)$ for all $a \in \mathbb{K}^{\times}$ (see [16]) is a formal analog of the connected component of the identity in a Lie group. We show that it inherits a structure of multipliability group from a group of the form $1 + \mathfrak{m}$.

**Theorem 2.** [Corollary 5.4] *The group $1\text{-}\mathrm{Aut}_k^+(\mathbb{K})$ has a natural structure of multipliability group.*

Section 5 is dedicated to the introduction of the formalism of "chain rules" on algebras of Noetherian series [11, 23, 13, 4]. Let $(\mathbb{K}, v)$ be an ordered valued field of Hahn series equipped with a derivation $\partial : \mathbb{K} \longrightarrow \mathbb{K}\,;\, a \mapsto a'$ and a composition law $\circ : \mathbb{K} \times \mathbb{K}^{>k} \longrightarrow \mathbb{K}$, and assume that $(\mathbb{K}^{>k}, \circ)$ forms a group with identity $x$. What then distinguishes operators $\circ_b : \mathbb{K} \longrightarrow \mathbb{K}$ ; $a \mapsto a \circ b$ from general field automorphisms of $\mathbb{K}$ is that the former satisfy the chain rule $(a \circ b)' = b' \cdot (a' \circ b)$. So $\circ_b$ is among automorphisms $\sigma$ for which there exists a multiplication operator $\mu : a \mapsto b' \cdot a$ with

$$\partial \circ \sigma = \mu \circ \sigma \circ \partial. \tag{2}$$

Using formal methods, we show that certain groups of automorphisms which satisfy a chain rule (2) are closed under infinite products.

**Theorem 3.** [Theorem 5.6] *Let $\phi \in \mathrm{Lin}^+_{\prec}(\mathbb{K})$ and assume that there is a $\xi \in \mathbb{K}$ with $\phi(\xi) = 1$ and $v(\phi(a)) > 0$ whenever $v(a) > v(\xi)$. Then the subgroup of $(1\text{-}\mathrm{Aut}^+_k(\mathbb{K}), \circ, \mathrm{Id}_{\mathbb{K}})$ of automorphisms $\sigma$ that satisfy a chain rule is closed under infinite products in $1\text{-}\mathrm{Aut}^+_k(\mathbb{K})$.*

In fact, the previous results are established in the more general setting of algebras of Noetherian series. We then apply this in Section 6 to the case when $\mathbb{K}$ is the field $\mathbb{S}$ of flat transseries. Flat transseries form a subfield of log-exp transseries as per [5, 7, 6] whose transmonomials are “flater than the identity”, i.e. are bounded by integer powers of the identity series $x$. Our main results are as follows:

**Theorem 4.** [Corollary 6.5] *The group $\mathcal{T} = (\{x + \delta \in \mathbb{S} : \mathbb{N}\,|\delta| < x\}, \circ, x)$ has a natural multipliability structure.*

For each $s \in \mathcal{T} \setminus \{x\}$, there is a natural isomorphism $r \mapsto s^{[r]}$ between $(k, +)$ and the centraliser $\mathcal{C}(s)$ of $s$ in $(\mathcal{T}, \circ)$. Given $s, t \in \mathcal{T}$, write $s \ll t$ if $\mathcal{C}(s) < \max(t, t^{\mathrm{inv}})$. The *residue* of an $s \in \mathcal{T} \setminus \{x\}$ corresponds to its equivalence class for the relation $s \sim t \Longleftrightarrow s \circ t^{\mathrm{inv}} \ll s$. Let $\mathcal{S}$ be a set of representatives for each residue in $\mathcal{T}$. Under a summability condition on $\mathcal{S}$, each $s \in \mathcal{T}$ can be represented uniquely as a well-ordered composition

$$s = \cdots \circ s_\gamma \circ \cdots \circ s_1 \circ s_0, \tag{3}$$

where $\lambda$ is an ordinal and $(s_\gamma)_{\gamma < \lambda}$ is a strictly decreasing sequence of elements of $\mathcal{S}$ for $\ll$:

**Theorem 5.** [Theorem 6.12] *There is a subset $S \subseteq \mathcal{T}$ such that for each $s \in \mathcal{T}$, there are a unique ordinal $\lambda$, a unique strictly $\ll$-decreasing function $s : \lambda \longrightarrow S$ and a unique function $r : \lambda \longrightarrow \mathbb{R}^\times$ with*

$$s = \cdots \circ s_\gamma^{[r(\gamma)]} \circ \cdots \circ s_1^{[r(1)]} \circ s_0^{[r(0)]}.$$

The choice of left-trailing well-ordered expansion (3) is somewhat arbitrary. We introduce a class of initial segments of linear orderings on ordinals (see Remark 1.3) along which transfinite compositions of elements of $\mathcal{S}$ are in one-to-one correspondence with elements of $\mathcal{T}$. The methods in Sections 5 and 6 are designed to adapt to larger fields than that of flat transseries (see Remark 6.13).

# 1 Groups and linearly ordered sets

## 1.1 Cartesian powers of groups

If $(\mathcal{G}, \cdot, 1)$ is a group and $X$ is a set, then we write $\mathcal{G}^X$ for the group, under pointwise product, of functions $I \longrightarrow \mathcal{G}$. That is, for all $f, g \in \mathcal{G}^X$ and $x \in X$, we have

$$(f \cdot g)(x) = f(x) \cdot g(x) \qquad \text{and} \qquad f^{-1}(x) = f(x)^{-1}.$$

The *support* of an $f\in\mathcal{G}^X$ is the set

$$\operatorname{supp} f := \{x\in X : f(x)\neq 1\}.$$

We denote by $\mathcal{G}^{(X)}$ the subset of elements $f\in\mathcal{G}^X$ whose support is finite, which is a normal subgroup of $\mathcal{G}^X$.

## 1.2 Linearly ordered sets

Recall that a linearly ordered set is an ordered pair $I=(\underline{I},<)$ where $\underline{I}$ is a set and $<$ is a binary relation on $\underline{I}$ such that for all $i,j,l\in\underline{I}$, we have

$$i\not< i \qquad \text{and} \qquad i<j\wedge j<l\Longrightarrow i<l \qquad \text{and} \qquad (i=j)\vee(i<j)\vee(j<i).$$

The category of all linearly ordered sets with nondecreasing maps as arrows is denoted by **Los**. *Isomorphisms* are to be understood in this sense, i.e. they are strictly increasing and bijective maps. Working in $I$, given an $i\in\underline{I}$, we write

$$\begin{aligned} i^- &:= (\{j\in\underline{I}: j<i\},<) \quad \text{and}\\ i^+ &:= (\{j\in\underline{I}: i<j\},<), \end{aligned}$$

expecting that the ambient linearly ordered set $I$ is clear from the context.

We now give a few conventions and notations pertaining to linear orderings.

- If $\mathcal{I}=(I,<)$ is a linearly ordered set, then we write $\underline{\mathcal{I}}:=I$ for its underlying set. We also write $\mathcal{I}^*=(I,<^*)$ where $<^*$ is the *reverse ordering*

  $$\forall i,j\in I,(i<^* j\Longleftrightarrow j<i).$$

  We also sometimes denote $<^*$ by $>$. Note that $\mathcal{I}^*\in\mathbf{Los}$ and that $(\mathcal{I}^*)^*=\mathcal{I}$.

- Let $I=(\underline{I},<)\in\mathbf{Los}$ and let $X\subseteq\underline{I}$ be a subset. We say that $X$ is a *convex subset* of $I$ if

  $$\forall i,j,l\in\underline{I},((i,l\in X\wedge i<j<l)\Longrightarrow j\in X).$$

  We say that $X$ is an *initial subset* of $I$ if

  $$\forall i,j\in\underline{I},((j\in X\wedge i<j)\Longrightarrow i\in X).$$

  We say that $X$ is a *final subset* of $I$ if it is an initial subset of $I^*$, i.e. if

  $$\forall i,j\in\underline{I},((j\in X\wedge j<i)\Longrightarrow i\in X).$$

- If $I=(\underline{I},<_I)$, $J=(\underline{J},<_J)$ are linearly ordered sets, then we write $I+J$ for the linearly ordered set $(\underline{I}\times\{0\}\sqcup\underline{J}\times\{1\},<)$ where $(i,0)<(i',0)\Longleftrightarrow i<_I i'$ and $(j,1)<(j',1)\Longleftrightarrow j<_J j'$ and $(i,0)<(j,1)$ for all $i\in\underline{I}$ and $j\in\underline{J}$. Note that the operation $+$ is associative up to isomorphism, i.e. there is a natural isomorphism between the functors $(\cdot+\cdot)+\cdot$ and $\cdot+(\cdot+\cdot)$ from $\mathbf{Los}^3$ to $\mathbf{Los}$. Therefore, for $n\in\mathbb{N}$ and $I_1,\ldots,I_n\in\mathbf{Los}$, we simply write $I_1+\cdots+I_n$ for $(\cdots((I_1+I_2)+I_3)+\cdots+I_{n-1})+I_n$.

- Let $J=(\underline{J},<)\in\mathbf{Los}$. Consider a linearly ordered set $I$ with ordering $<_I$ and a family $(J_i)_{i\in\underline{I}}$ of linearly ordered subsets $(\underline{J_i},<)$ of $J$. We write

  $$J=\coprod_I (J_i)_{i\in\underline{I}}$$

if we have

$$\underline{J} = \bigcup_{i \in I} \underline{J_i} \qquad \text{and} \qquad \forall i, i' \in I, (i <_I i' \Longrightarrow (\forall (j, j') \in \underline{J_i} \times \underline{J_{i'}}, (j < j'))).$$

Then $\underline{J}$ is the disjoint union $\underline{J} = \bigsqcup_{i \in I} \underline{J_i}$. Any choice function in $\prod_{i \in \underline{I}} \underline{J_i}$ is a strictly increasing map $I \longrightarrow J$. Conversely, if $L$ is a subset of $\underline{J}$, then writing

$$\underline{J_l} := \{j \in \underline{J} : ((j, l) \cup (l, j)) \cap L = \varnothing\}$$

for each $l \in L$, each $\underline{J_l}$ is a convex subset of $J$ and we have $J = \coprod_{(L,<)} (J_l)_{l \in L}$.

## 1.3 Tree-like orderings on ordinals

We fix a non-zero ordinal $\lambda$, write $\lambda - 1 := \{\beta \in \lambda : \beta + 1 < \lambda\}$, so $\lambda - 1 = \lambda$ if $\lambda$ is a limit and $\lambda - 1 = \beta$ if $\lambda = \beta + 1$ is a successor. We will define a specific linear ordering on $\lambda$ that is suited to perform inductive definitions of products over ill-ordered linear orderings.

Given a map $N : \lambda - 1 \longrightarrow \{-1, 1\}$, we define a binary relation $<_N$ on $\lambda$ as follows. For $\alpha, \beta \in \lambda$ with $\alpha < \beta$, we set $\alpha <_N \beta$ if $N(\alpha) = 1$, and we set $\beta <_N \alpha$ if $N(\alpha) = -1$. In other words, we have $\alpha <_N \beta$ if and only if

$$(N \upharpoonright \alpha) \cup ([\alpha, \lambda - 1) \times \{0\}) <_{\text{lex}} (N \upharpoonright \beta) \cup ([\beta, \lambda - 1) \times \{0\})$$

for the lexicographic ordering on functions $\lambda - 1 \longrightarrow \{-1, 0, 1\}$. It follows that $<_N$ is a linear ordering on $\lambda$.

**Lemma 1.1.** *For each $\alpha < \mu < \lambda$, the set $[\alpha, \mu)$ is a convex subset of $(\mu, <_N)$.*

**Proof.** By the characterisation of $<_N$ by $<_{\text{lex}}$ above, we have

$$[\alpha, \mu) = \{\beta \in \mu : \forall \gamma < \alpha, ((\gamma <_N \alpha \Longrightarrow \gamma <_N \beta) \wedge (\alpha <_N \gamma \Longrightarrow \beta <_N \alpha))\},$$

which is clearly a convex subset of $(\mu, <_N)$. □

For $\alpha < \mu < \lambda$, we write

$$L_N(\alpha, \mu) := \{\beta < \mu : \beta <_N [\alpha, \mu)\} \qquad \text{and} \qquad R_N(\alpha, \mu) := \{\beta < \mu : [\alpha, \mu) <_N \beta\},$$

so $(\mu, <_N) = (L_N(\alpha, \mu), <_N) \amalg ([\alpha, \mu), <_N) \amalg (R_N(\alpha, \mu), <_N)$.

**Lemma 1.2.** *For $\beta < \alpha$, the set $L_N(\beta, \mu)$ is an initial segment of $(L_N(\alpha, \mu), <_N)$ whereas $R_N(\beta, \mu)$ is a final segment of $(R_N(\alpha, \mu), <_N)$.*

**Proof.** Let $\gamma \in L_N(\beta, \mu)$ and $\eta \in L_N(\alpha, \mu)$ with $\eta <_N \gamma$. Since $\gamma \in L_N(\beta, \mu)$, we have $\gamma <_N \beta$ so $\eta <_N \beta$, whence in particular $\eta \notin R_N(\beta, \mu)$. Assume for contradiction that $\eta \notin L_N(\beta, \mu)$. So $\eta \in [\beta, \mu)$. But $[\beta, \mu) \supseteq [\alpha, \mu)$ where $L_N(\alpha, \mu) \cap [\alpha, \mu) = \varnothing$: a contradiction. So $\eta \in L_N(\beta, \mu)$ which is thus an initial segment of $(L_N(\alpha, \mu), <_N)$. The proof that $R_N(\beta, \mu)$ is a final segment of $(R_N(\alpha, \mu), <_N)$ is symmetric. □

**Remark 1.3.** Conversely, if $\lhd$ is a linear ordering on $\lambda$ satisfying the conclusion of Lemma 1.1, then $\lhd$ coincides with $<_N$ for the function $N : \lambda - 1 \longrightarrow \{-1, 1\}$ defined by $N(\alpha) = 1$ if $\alpha \lhd [\alpha + 1, \lambda)$ and $N(\alpha) = -1$ if $[\alpha + 1, \lambda) \lhd \alpha$.

# 2 Multipliability groups

## 2.1 Groups with linearly ordered products

Let $(\mathcal{G},\cdot,1)$ be a group. We wish to extend the elementary properties of $\Pi^{\min}$ to families $\Pi=(\Pi_I)_{I\in\mathbf{Los}}$ whose domains $\operatorname{dom}\Pi_I$ may be strictly larger than $\mathcal{G}^{(I)}$ for certain $I\in\mathbf{Los}$. This leads to the following formalism. Consider a family $\Pi=(\Pi_I)_{I\in\mathbf{Los}}$ of partial functions

$$\Pi_I:\mathcal{G}^{\underline{I}}\longrightarrow\mathcal{G}$$

whose domains are denoted $\operatorname{dom}\Pi_I$. Given $I,J\in\mathbf{Los}$, we regard the following axioms:

**MG1.** The set $\operatorname{dom}\Pi_I$ is a subgroup of $\mathcal{G}^{\underline{I}}$. Moreover, for $f,g\in\operatorname{dom}\Pi_I$, writing

$$f[g]:=(\Pi_{(i^+)^*}(f^{-1}\upharpoonright(\underline{i^+}))\cdot g(i)\cdot\Pi_{i^+}(f\upharpoonright(\underline{i^+})))_{i\in\underline{I}},$$

we have

$$\Pi_I(f\cdot g)=(\Pi_I f)\cdot(\Pi_I f[g]).$$

**MG2.** We have $\operatorname{dom}\Pi_I\supseteq\mathcal{G}^{(I)}$, and for all $g\in\mathcal{G}^{(I)}$, writing $\operatorname{supp} g=\{i_1,\ldots,i_n\}$ where $n\in\mathbb{N}$ and $i_1<\cdots<i_n$, we have

$$\Pi_I g=g(i_1)\cdots g(i_n).$$

**MG3.** If $\varphi:I\longrightarrow J$ is an isomorphism and $g\in\operatorname{dom}\Pi_J$, then $g\circ\varphi\in\operatorname{dom}\Pi_I$, and

$$\Pi_I(g\circ\varphi)=\Pi_J g.$$

**MG4.** Let $(I_j)_{j\in J}$ be a family with $I=\coprod_J(I_j)_{j\in\underline{J}}$. Let $g\in\operatorname{dom}\Pi_I$ and for each $j\in\underline{J}$, set $g_j:=g\upharpoonright\underline{I_j}$. Then

**MG4a.** $g_j\in\operatorname{dom}\Pi_{I_j}$ for each $j\in\underline{J}$.

**MG4b.** $(\Pi_{I_j}g_j)_{j\in\underline{J}}\in\operatorname{dom}\Pi_J$.

**MG4c.** $\Pi_J((\Pi_{I_j}g_j)_{j\in\underline{J}})=\Pi_I g$.

**MG5.** If $I=I_1\amalg I_2$ and $(g,h)\in\operatorname{dom}\Pi_{I_1}\times\operatorname{dom}\Pi_{I_2}$, then $(g\sqcup h)\in\operatorname{dom}\Pi_I$.

**MG6.** Let $g\in\operatorname{dom}\Pi_I$ and $g_0\in\mathcal{G}$. Then $g_0\cdot g\cdot g_0^{-1}=(g_0\cdot g(i)\cdot g_0^{-1})_{i\in\underline{I}}\in\operatorname{dom}\Pi_I$ and

$$\Pi_I(g_0\cdot g\cdot g_0^{-1})=g_0\cdot(\Pi_I g)\cdot g_0^{-1}.$$

**MG7.** Let $g\in\operatorname{dom}\Pi_I$. Then $g^{-1}\in\operatorname{dom}\Pi_{I^*}$ and

$$\Pi_{I^*}g^{-1}=(\Pi_I g)^{-1}.$$

**Fullness.** For each linear ordering $\lhd$ on $\underline{I}$, we have $\operatorname{dom}\Pi_I=\operatorname{dom}\Pi_{(\underline{I},\lhd)}$.

**Definition 2.1.** *If $(\mathcal{G},\cdot,1)$ is a group and $\Pi$ is a function satisfying* **MG1–MG7**, *then we say that $\Pi$ is* **multipliability structure** *on $(\mathcal{G},\cdot,1)$ and that $(\mathcal{G},\Pi)$ is a* **multipliability group**. *If furthermore $\Pi$ satisfies* **Fullness**, *then we say that $(\mathcal{G},\Pi)$ (resp. $\Pi$) is* **full**.

**Remark 2.2.** If $(\mathcal{G},\cdot,1,\Pi)$ satisfies **MG4**, then for $I=(\underline{I},<)\in\mathbf{Los}$, $f\in\operatorname{dom}\Pi_I$ and $X\subseteq\underline{I}$ such that $(X,<)$ is a convex subset of $I$, then we have $f\upharpoonright X\in\operatorname{dom}\Pi_{(X,<)}$. Indeed, we can write $I=(L,<)\amalg(X,<)\amalg(R,<)$ where

$$(L,R):=(\{i\in\underline{I}:\forall x\in X,i<X\},\{i\in\underline{I}:\forall x\in X,i>x\}).$$

Therefore, this follows from **MG4a**.

**Remark 2.3.** In **MG1**, it is implicitly assumed that $\Pi_{(i^+)^*}(f^{-1} \upharpoonright (\underline{i^+}))$, $\Pi_{i^+}(f \upharpoonright (\underline{i^+}))$ and $\Pi_I f[g]$ are well-defined. This follows from Remark 2.2 once one knows that $(\mathcal{G}, \cdot, 1, \Pi)$ satisfies **MG4** and **MG7**. So when showing that $\Pi$ is a multipliability structure, we will omit the proofs of well-definedness of those families, and directly prove the other axioms.

We say that a family $g : \underline{I} \longrightarrow \mathcal{G}$ is *$I$-multipliable* if $g \in \operatorname{dom} \Pi_I$. If $\Pi$ is full, then we simply say that $g$ is *multipliable* if it is $I$-multipliable.

**Definition 2.4.** *Let $(\mathcal{G}, \cdot, 1, \Pi^{\mathcal{G}})$ and $(\mathcal{H}, \cdot, 1, \Pi^{\mathcal{H}})$ be strong groups. A* **strongly multiplicative group morphism** *is a group morphism $\Phi : \mathcal{G} \longrightarrow \mathcal{H}$ such that for for all $I \in \mathbf{Los}$ and $g \in \operatorname{dom} \Pi_I^{\mathcal{G}}$, we have $\Phi \circ g \in \operatorname{dom} \Pi_I^{\mathcal{H}}$ and $\Pi_I^{\mathcal{H}}(\Phi \circ g) = \Phi(\Pi_I^{\mathcal{G}} g)$.*

We write $\operatorname{Hom}^+(\mathcal{G}, \mathcal{H})$ for the set of strongly multiplicative group morphisms $\mathcal{G} \longrightarrow \mathcal{H}$ and $\operatorname{Aut}^+(\mathcal{G})$ for the set of strongly multiplicative group isomorphisms $\mathcal{G} \longrightarrow \mathcal{G}$. Note that the compositum of composable strongly multiplicative group morphisms is strongly multiplicative. In particular $\operatorname{Aut}^+(\mathcal{G})$ is a subgroup of the group of automorphisms of $\mathcal{G}$.

## 2.2 Examples

Throughout Section 2.2, we fix a group $(\mathcal{G}, \cdot, 1)$.

**Example 2.5.** It is routine to check that $\Pi^{\min}$ is a full multipliability structure on $\mathcal{G}$. We call $\Pi^{\min}$ the *minimal multipliability structure* on $\mathcal{G}$. Note that if $(\mathcal{H}, \cdot, 1, \Pi)$ is a multipliability group and $\mathcal{G}$ is equipped with $\Pi^{\min}$, then $\operatorname{Hom}^+(\mathcal{G}, \mathcal{H}) = \operatorname{Hom}(\mathcal{G}, \mathcal{H})$.

**Definition 2.6.** *Let $(\mathcal{G}, \cdot, 1, \Pi)$ be a multipliability group. A* **multipliability subgroup** *of $(\mathcal{G}, \Pi)$ is a subgroup $\mathcal{H} \subseteq \mathcal{G}$ such that $\Pi_I h \in \mathcal{H}$ whenever $I \in \mathbf{Los}$ and $h \in \operatorname{dom} \Pi_I \cap \mathcal{H}^{\underline{I}}$.*

**Proposition 2.7.** *Let $(\mathcal{G}, \cdot, 1, \Pi)$ be a multipliability group and let $\mathcal{H}$ be a multipliability subgroup. The function $\Pi^{\mathcal{H}}$ given by $\Pi_I^{\mathcal{H}} := \Pi_I \upharpoonright \mathcal{H}^{\underline{I}}$ for each $I \in \mathbf{Los}$ is a multipliability structure on $\mathcal{H}$. Moreover, it is full if $\Pi$ is full.*

**Proof.** It is easy to see that the axioms are preserved. □

Thus multipliability subgroups are the right notion of substructure in our framework.

## 2.3 Quotients

Let $(\mathcal{G}, \cdot, 1, \Pi)$ be a multipliability group and let $(\mathcal{H}, \cdot, 1)$ be a group. Let $\Phi : \mathcal{G} \longrightarrow \mathcal{H}$ be a group morphism and suppose that $\operatorname{Ker}(\Phi)$ is a multipliability subgroup of $(\mathcal{G}, \Pi)$. We define a multipliability structure $\Pi^{\Phi}$ on $\Phi(\mathcal{G}) \subseteq \mathcal{H}$ as follows. For all linearly ordered sets $I$, define

$$\operatorname{dom} \Pi_I^{\Phi} = \{\Phi \circ g : g \in \operatorname{dom} \Pi_I\} \qquad \text{and} \qquad \Pi_I^{\Phi}(\Phi \circ g) := \Phi(\Pi_I g)$$

for all $g \in \operatorname{dom} \Pi_I$. If $g' \in \operatorname{dom} \Pi_I$ is such that $\Phi \circ g = \Phi \circ g'$, then we have $g^{-1} \cdot g' \in \operatorname{Ker}(\Phi)^{\underline{I}}$. Now $g' = g \cdot (g^{-1} \cdot g')$ so **MG1** gives $\Pi_I g' = (\Pi_I g) \cdot (\Pi_I g[g^{-1} \cdot g'])$. Since $\operatorname{Ker}(\Phi)$ is a normal subgroup of $\mathcal{G}$, the family $g[g^{-1} \cdot g']$ lies in $\operatorname{Ker}(\Phi)^{\underline{I}}$. It follows since $\operatorname{Ker}(\Phi)$ is a multipliability subgroup of $(\mathcal{G}, \Pi)$ that $\Pi_I g[g^{-1} \cdot g'] \in \operatorname{Ker}(\Phi)$, so $\Phi(\Pi_I g') = \Phi(\Pi_I g) \cdot \Phi(\Pi_I g[g^{-1} \cdot g']) = \Phi(\Pi_I g)$. Therefore $\Pi_I^{\Phi}$ is well-defined.

**Proposition 2.8.** *The structure* $(\Phi(\mathcal{G}), \Pi^{\Phi})$ *is a multipliability group. It is full if* $\Pi$ *is full.*

**Proof.** We fix linearly ordered sets $I, J$. Let $q = \Phi \circ f, r = \Phi \circ g \in \operatorname{dom} \Pi_I^{\Phi}$ where $f, g \in \operatorname{dom} \Pi_I$. Then $q \cdot r^{-1} = \Phi \circ (f \cdot g^{-1})$ where $f \cdot g^{-1} \in \operatorname{dom} \Pi_I$, so $q \cdot r^{-1} \in \operatorname{dom} \Pi_I^{\Phi}$. So $\operatorname{dom} \Pi_I^{\Phi}$ is a subgroup of $\Phi(\mathcal{G})^{\underline{I}}$. For $s = \Phi \circ h \in \Phi(\mathcal{G})^{\underline{I}}$ where $h \in \mathcal{G}^{\underline{I}}$, we have $s \cdot q \cdot s^{-1} = \Phi \circ (h \cdot f \cdot h^{-1})$ where $h \cdot f \cdot h^{-1} \in \operatorname{dom} \Pi_I$ by **MG1**. So $s \cdot q \cdot s^{-1} \in \operatorname{dom} \Pi_I^{\Phi}$, whence $\operatorname{dom} \Pi_I^{\Phi}$ is a normal subgroup of $\Phi(\mathcal{G})^{\underline{I}}$. For $i \in \underline{I}$, we have

$$\begin{aligned} q[r](i) &= \Pi_{(i^+)^*}^{\Phi}(\Phi \circ f^{-1} \upharpoonright (\underline{i^+})) \cdot (\Phi \circ g)(i) \cdot \Pi_{i^+}^{\Phi}(\Phi \circ f \upharpoonright (\underline{i^+})) \\ &= \Phi(\Pi_{(i^+)^*}(f^{-1} \upharpoonright (\underline{i^+}))) \cdot \Phi(g(i)) \cdot \Phi(\Pi_{i^+}(f \upharpoonright (\underline{i^+}))) \\ &= \Phi(f[g](i)). \end{aligned}$$

Therefore **MG1** in $\mathcal{G}$ gives

$$(\Pi_I^{\Phi} q) \cdot (\Pi_I^{\Phi} q[r]) = \Phi(\Pi_I q) \cdot \Phi(\Pi_I f[g]) = \Phi((\Pi_I f) \cdot (\Pi_I f[g])) = \Phi(\Pi_I f \cdot g) = \Pi_I^{\Phi}(q \cdot r).$$

The axioms **MG2**, **MG3** and **MG5**–**MG7** follow immediately from the validity of their respective counterparts in $(\mathcal{G}, \Pi)$.

Assume that $I = \coprod_J (I_j)_{j \in \underline{J}}$ and write $q_j := q \upharpoonright \underline{I_j}$ and $f_j := f \upharpoonright \underline{I_j}$ for all $j \in \underline{J}$. Then by **MG4** in $(\mathcal{G}, \Pi)$, we have that

- for $j \in \underline{J}$, $f_j$ lies in $\operatorname{dom} \Pi_{I_j}$, so $q_j$ lies in $\operatorname{dom} \Pi_{I_j}^{\Phi}$, and we have

$$\Pi_{I_j}^{\Phi} q_j = \Phi(\Pi_{I_j} f_j). \tag{2.1}$$

- the family $(\Pi_{I_j} f_j)_{j \in \underline{J}}$ lies in $\operatorname{dom} \Pi_J$, so $(\Pi_{I_j}^{\Phi} q_j)_{j \in \underline{J}} \in \operatorname{dom} \Pi_J^{\Phi}$ by (2.1).
- $\Pi_J (\Pi_{I_j} f_j)_{j \in \underline{J}} = \Pi_I f$, so

$$\begin{aligned} \Pi_J^{\Phi} (\Pi_{I_j}^{\Phi} q_j)_{j \in \underline{J}} &= \Pi_J^{\Phi} (\Phi(\Pi_{I_j} f_j))_{j \in \underline{J}} \\ &= \Phi(\Pi_J (\Pi_{I_j} f_j)_{j \in \underline{J}}) \\ &= \Phi(\Pi_I f) \\ &= \Pi_I^{\Phi} q. \end{aligned}$$

This proves **MG4**. If $(\mathcal{G}, \Pi)$ is full, then one easily sees that $(\Phi(\mathcal{G}), \Pi^{\Phi})$ is full as well. □

Note that $\Phi : (\mathcal{G}, \Pi) \longrightarrow (\Phi(\mathcal{G}), \Pi^{\Phi})$ is strongly multiplicative by definition. Applying this for quotient maps and injective morphisms respectively, we obtain:

**Corollary 2.9.** *If* $\mathcal{H}$ *is a normal multipliability subgroup of* $\mathcal{G}$*, then* $(\mathcal{G}/\mathcal{H}, \Pi^{\Phi})$ *is a multipliability group for the quotient map* $\Phi : \mathcal{G} \twoheadrightarrow \mathcal{G}/\mathcal{H}$*. Moreover* $(\mathcal{G}/\mathcal{H}, \Pi^{\Phi})$ *is full if* $(\mathcal{G}, \Pi)$ *is full.*

**Corollary 2.10.** *If* $\Psi : \mathcal{G} \longrightarrow \mathcal{H}$ *is an injective morphism, then* $(\Psi(\mathcal{G}), \Pi^{\Psi})$ *is a multipliability group which is full if* $(\mathcal{G}, \Pi)$ *is full.*

# 3 Strongly linear algebra

In this section, we recall facts regarding a commutative version of multipliability groups that was introduced in [12, Section 6.2] and expanded upon in [23, 8, 4]. We fix a field $k$.

## 3.1 Summability spaces and strongly linear functions

Let $(V, +, 0, .)$ be a vector space over $k$. Let us be given, for each set $I$, a $k$-linear function

$$\Sigma_I : \mathrm{dom}\,\Sigma_I \longrightarrow V$$

whose domain $\mathrm{dom}\,\Sigma_I$ is a subspace of the vector space $V^I$. Given sets $I, J$, consider the following axioms:

**SS1.** $V^{(I)} \subseteq \mathrm{dom}\,\Sigma_I$ and $\Sigma_I v = \sum_{i \in \mathrm{supp}\, v} v(i)$ for all $v \in V^{(I)}$.

**SS2.** if $\varphi : I \longrightarrow J$ is bijective and $v \in \mathrm{dom}\,\Sigma_J$, then $v \circ \varphi \in \mathrm{dom}\,\Sigma_I$, and $\Sigma_I(v \circ \varphi) = \Sigma_J v$.

**SS3.** if $I = \bigsqcup_{j \in J} I_j$ and $v \in \mathrm{dom}\,\Sigma_J$, then writing $v_j := v \upharpoonright I_j$ for each $j \in J$, we have

**SS3a.** $v_j \in \mathrm{dom}\,\Sigma_{I_j}$ for all $j \in J$,

**SS3b.** $\left(\sum_{I_j} v_j\right)_{j \in J} \in \mathrm{dom}\,\Sigma_J$, and

**SS3c.** $\sum_I v = \sum_J \left(\left(\sum_{I_j} v_j\right)_{j \in J}\right)$.

**SS4.** if $I = I_1 \sqcup I_2$ and $(v, w) \in \mathrm{dom}\,\Sigma_{I_1} \times \mathrm{dom}\,\Sigma_{I_2}$, then the family $(v \sqcup w) : I \longrightarrow V$ given by $(v \sqcup w)(i_1) := v(i_1)$ and $(v \amalg w)(i_2) = w(i_2)$ for all $i_1 \in I_1$ and $i_2 \in I_2$ lies in $\mathrm{dom}\,\Sigma_I$.

**UF.** For all $v \in \mathrm{dom}\,\Sigma_I$ and all families $(f_i)_{i \in I}$ of scalar valued functions $f_i : X_i \longrightarrow k$ with finite domains $X_i$, writing $I' = \{(i, x) : i \in I \wedge x \in X_i\}$, the family $(f_i(x)\, v(i))_{i \in I \wedge x \in X_i}$ lies in $\mathrm{dom}\,\Sigma_{I'}$.

We say that $\Sigma = (\Sigma_I)_{I \in \mathbf{Set}}$ is a *summability structure* on $V$, or that $(V, \Sigma)$ is a *summability space*, if $(V, \Sigma)$ satisfies **SS1**–**SS4** and **UF** for all sets $I$.

**Remark 3.1.** Here, we call summability spaces what are called ultrafinite summability spaces in [4].

Given a set $I$, a family $v \in V^I$ is said *summable* (in $(V, \Sigma)$) if $v \in \mathrm{dom}\,\Sigma_I$. A subspace $W \subseteq V$ is a *summability subspace* if for all sets $I$ and families $v \in \mathrm{dom}\,\Sigma_I \cap W^I$, we have $\Sigma_I v \in W$. Then we have an induced summability structure $\Sigma^W$ on $W$ defined by $\Sigma_I^W := \Sigma_I \cap (W^I \times W)$ for each set $I$.

Let $(V, \Sigma)$ and $(V_1, \Sigma_1)$ be summability spaces over $k$. A function $\phi : V \longrightarrow V_1$ is said *strongly linear* if it is linear, and if, for all sets $I$ and families $v \in \mathrm{dom}\,\Sigma_I$, the family $\phi \circ v$ lies in $\mathrm{dom}\,\Sigma_{1,I}$ and $\Sigma_{1,I}(\phi \circ v) = \phi(\Sigma_I v)$. We write $\mathrm{Lin}^+(V)$ for the vector space of strongly linear functions $V \longrightarrow V$.

A family $\phi : J \longrightarrow \mathrm{Lin}^+(V)$ is said Lin-*summable* if for all sets $I$ and all families $v \in \mathrm{dom}\,\Sigma_I$, we have $\phi(v) := (\phi(j)(v(i)))_{(i,j) \in I \times J} \in \mathrm{dom}\,\Sigma_{I \times J}$. For $v_0 \in V$, the family $(v_0)_{i \in \{0\}}$ is summable by **SS1**, so $\phi(v_0) := (\phi(j)(v_0))_{j \in J}$ is summable, and we set

$$\left(\sum_{j \in J} \phi(j)\right)(v_0) := \Sigma_J\, \phi(v_0).$$

One can check that the function $\sum_{j \in J} \phi(j)$ thus defined is strongly linear. Lin-summability defines [4, Proposition 1.25] a summability structure on $\mathrm{Lin}^+(V)$.

## 3.2 Summability algebras

Let $(A, +, \cdot, .)$ be a $k$-algebra, where $(A, +, \cdot)$ is a possibly non-commutative, possibly non-unital ring. We consider $A$ as a Lie algebra under the Lie bracket $[a, b] := a \cdot b - b \cdot a$, and understand Lie subalgebras of $A$ in this sense.

Let $\Sigma$ be a summability structure on $(A, +, .)$. We say that $(A, \Sigma)$ is a *summability algebra* if the following axiom is satisfied.

**SA.** For all sets $I$, $J$ and all $(a, b) \in \operatorname{dom} \Sigma_I \times \operatorname{dom} \Sigma_J$, the family $a \cdot b := (a(i) \cdot b(j))_{(i,j) \in I \times J}$ lies in $\operatorname{dom} \Sigma_{I \times J}$, and we have

$$\sum_{I \times J} (a \cdot b) = \left(\sum_I a\right) \cdot \left(\sum_J b\right).$$

Given a summability space $(V, \Sigma)$, the algebra $(\operatorname{Lin}^+(V), +, ., \circ)$, with its summability structure of Lin-summability, is a summability algebra [4, Proposition 1.30].

**Proposition 3.2.** [4, Proposition 1.28] *Let $(A, \Sigma)$ be a summability algebra. Then for all $a \in A$, the left and right product functions $a \cdot : A \longrightarrow A; b \mapsto a \cdot b$ and $\cdot a : A \longrightarrow A; b \mapsto b \cdot a$ are strongly linear.*

Let $(A, \Sigma)$ be a summability $k$-algebra. A *summability ideal* of a summability algebra $(A, \Sigma)$ is a *two-sided* ideal of $A$ which is a summability subspace of $(A, \Sigma)$. A s*ummability Lie subalgebra* of $\operatorname{Lin}^+(A)$ is a Lie subalgebra of $\operatorname{Lin}^+(A)$, for the standard Lie bracket, that is also a summability subspace.

We write $\operatorname{End}^+(A)$ for the set of strongly linear endomorphisms of algebra of $A$. A *derivation* on $A$ is a $k$-linear function $\partial : A \longrightarrow A$ which satisfies the Leibniz product rule. We write $\operatorname{Der}^+(A)$ for the set of strongly linear derivations on $A$. It is a summability Lie subalgebra of $\operatorname{Lin}^+(A)$.

## 3.3 Formal power series

An important example of summability algebra is Magnus' algebra of formal series in non-commuting variables. Let $J$ be a set. Write

$$J^\star := \bigcup_{n \in \mathbb{N}} J^n$$

where $J^0 = \{\varnothing\}$. We see elements of $J^\star$ as *finite words* with letters in $J$. For $m, n \in \mathbb{N}$ and $(\beta, \gamma) = ((\beta_1, \ldots, \beta_m), (\gamma_1, \ldots, \gamma_n)) \in J^m \times J^n$, we define

$$\beta : \gamma := (\beta_1, \ldots, \beta_m, \gamma_1, \ldots, \gamma_n) \in J^{m+n},$$

where it is implied that $\varnothing : \theta = \theta : \varnothing = \theta$ for all $\theta \in J^\star$. Then $(J^\star, :, \varnothing)$ is a cancellative monoid.

We write $k\langle\langle J \rangle\rangle := k^{J^\star}$ for the vector space over $k$ of functions $J^\star \longrightarrow k$. This comes with a natural summability structure of pointwise summability, whereby a family $Q : I \longrightarrow k\langle\langle J \rangle\rangle$ is pointwise summable if for all $\theta \in J^\star$, the set $I_\theta := \{i \in I : \theta \in \operatorname{supp} Q_i\}$ is finite. We then define $\Sigma_I Q$ as the pointwise sum

$$\Sigma_I Q : J^\star \longrightarrow k; \theta \mapsto \sum_{i \in I_\theta} Q_i(\theta).$$

The space $k\langle\langle J \rangle\rangle$ is equipped with a Cauchy product: for $P, Q \in k\langle\langle J \rangle\rangle$, we have

$$\forall \theta \in J^\star, (P \cdot Q)(\theta) := \sum_{\theta = \beta : \gamma} P(\beta)\, Q(\gamma).$$

Note that

$$\begin{aligned} \operatorname{supp}(P + Q) &\subseteq \operatorname{supp} P \cup \operatorname{supp} Q \quad \text{and} \\ \operatorname{supp}(P \cdot Q) &\subseteq (\operatorname{supp} P) : (\operatorname{supp} Q). \end{aligned} \tag{3.1}$$

Then $k\langle\langle J\rangle\rangle$ is a unital summability algebra. Moreover the set

$$k\langle\langle J\rangle\rangle_0 := \{P \in k\langle\langle J\rangle\rangle : P(\varnothing) = 0\}$$

is a summability ideal.

For $\theta \in J^\star$, we write $X_\theta$ for the function $J^\star \longrightarrow k$ with support $\{\theta\}$ and $X_\theta(\theta) = 1$. So $X_\varnothing = 1$, and writing $\theta = (\theta_1, \ldots, \theta_n)$, we have $X_\theta = X_{\theta_1} \cdots X_{\theta_n}$. For all $P \in k\langle\langle J\rangle\rangle$, the family $(P(\theta)\, X_\theta)_{\theta \in J^\star}$ is pointwise summable, with $P = \sum_{\theta \in J^\star} P(\theta)\, X_\theta$. Lastly, for each subset $X \subseteq J$, the restriction map

$$\begin{aligned} \rho_X : k\langle\langle J\rangle\rangle &\longrightarrow k\langle\langle X\rangle\rangle \\ P &\longmapsto P \upharpoonright X^\star \end{aligned}$$

is a strongly linear morphism of algebras.

## 3.4 Evaluations

Let $(A, \Sigma)$ be a unital summability algebra of the form $A = k + \mathfrak{m}$ where $\mathfrak{m}$ is a summability ideal. Then $A$ *has evaluations* if for all sets $J$ and all $f \in \operatorname{dom} \Sigma_J^{\mathfrak{m}}$, the family $(f(\theta_1) \cdots f(\theta_n))_{\theta = (\theta_1, \ldots, \theta_n) \in J^\star}$ is summable. Then for all sets $J$, all $f \in \operatorname{dom} \Sigma_J^{\mathfrak{m}}$ and all $P \in k\langle\langle J\rangle\rangle$, we define the *evaluation* of $P$ at $f$ as

$$\operatorname{ev}_f(P) := \sum_{\theta = (\theta_1, \ldots, \theta_n) \in J^\star} P(\theta)\, f(\theta_1) \cdots f(\theta_n) \in A.$$

If $J = \{0, \ldots, m-1\}$ where $m \in \mathbb{N}$, then we simply write $\operatorname{ev}_{f(0), \ldots, f(m-1)}(P) := \operatorname{ev}_f(P)$. Each such evaluation map $\operatorname{ev}_f : k\langle\langle J\rangle\rangle \longrightarrow A$ is [4, Theorem 2.6] a strongly linear morphism of algebras. The summability algebra $k\langle\langle J\rangle\rangle = k + k\langle\langle J\rangle\rangle_0$ itself has evaluations.

**Proposition 3.3.** *Let $A = k + \mathfrak{m}$ be a summability algebra with evaluations. Then $A$ is local algebra with maximal ideal $\mathfrak{m}$. In particular, the set $(1 + \mathfrak{m}, \cdot, 1)$ is a subgroup of $A^\times$.*

The following proposition is very useful in order to derive identities in summability algebras from formal identities in $k\langle\langle J\rangle\rangle$.

**Proposition 3.4.** [4, Proposition 2.7] *Let $(A, \Sigma)$ be a summability algebra with evaluations, with maximal ideal $\mathfrak{m}$. Let $I, J$ be sets and let $Q : I \longrightarrow k\langle\langle J\rangle\rangle_0$ and $a : J \longrightarrow \mathfrak{m}$ be summable. For all $P \in k\langle\langle I\rangle\rangle$, we have*

$$\operatorname{ev}_a(\operatorname{ev}_Q(P)) = \operatorname{ev}_{(\operatorname{ev}_a(Q(i)))_{i \in I}}(P).$$

If $k$ has characteristic zero, then we have a formal power series

$$\log(1 + X_\bullet) := \sum_{m > 0} \frac{(-1)^{m-1}}{m} X_\bullet^m \in k\langle\langle\{\bullet\}\rangle\rangle_0.$$

**Proposition 3.5.** [4, Corollary 2.10] *Suppose that $k$ has characteristic zero. Let $(k + \mathfrak{m}, \Sigma)$ be a summability algebra with evaluations. Then we have a bijection*

$$\begin{aligned} \log : 1 + \mathfrak{m} &\longrightarrow \mathfrak{m} \\ \delta &\longmapsto \operatorname{ev}_{\delta - 1}(\log(1 + X_\bullet)) = \sum_{m > 0} \frac{(-1)^{m-1}}{m} (\delta - 1)^m. \end{aligned}$$

**Proposition 3.6.** *Suppose that $k$ has characteristic zero. Let $(A,\Sigma)$ be a summability algebra with evaluations, with maximal ideal $\mathfrak{m}$. Let $\varepsilon\in\mathfrak{m}$. The function*

$$\begin{array}{rcl} k & \longrightarrow & 1+\mathfrak{m}\\ c & \longmapsto & (1+\varepsilon)^{[c]} := \sum_{m\in\mathbb{N}} \binom{c}{m}\varepsilon^m \end{array}$$

*is a group morphism $(k,+)\longrightarrow(1+\mathfrak{m},\cdot)$ with $(1+\varepsilon)^{[1]}=1+\varepsilon$ and $((1+\varepsilon)^{[c]})^{[c']}=(1+\varepsilon)^{[cc']}$ for all $\varepsilon\in\mathfrak{m}$ and $c,c'\in k$.*

**Proof.** In view of [4, Proposition 2.16], it suffices to note that the reciprocal exp of log satisfies the classical formula $\exp(c\log(1+X_\bullet))=\sum_{m\in\mathbb{N}}\binom{c}{m}X_\bullet^m$ in $k\langle\langle\{\bullet\}\rangle\rangle$. □

**Remark 3.7.** This law turns $1+\mathfrak{m}$ into a $k$-group as per Miasnikov-Remeslennikov [21, 22].

# 4 Ordered products in algebras with evaluations

## 4.1 The formal ordered product

We fix two linearly ordered sets $J=(\underline{J},<)$ and $I=(\underline{I},<)$. Let $J^\star$ denote the subset of $(\underline{J})^\star$ of words $(\theta_1,\ldots,\theta_n)\in(\underline{J})^\star$ with $\theta_1<\cdots<\theta_n$. It is understood that $\varnothing\in J^\star$. Note that if $J=\coprod_I(J_i)_{i\in\underline{I}}$ for a family $(J_i)_{i\in\underline{I}}$, then we have

$$J^\star=\bigsqcup_{(i_1,\ldots,i_n)\in I^\star}(J_{i_1})^\star:\cdots:(J_{i_n})^\star. \tag{4.1}$$

Define

$$\mathrm{OP}_J:=\sum_{\theta\in J^\star}X_\theta\in k\langle\langle\underline{J}\rangle\rangle.$$

Note that $\varnothing^\star=\{\varnothing\}$ and that $\mathrm{OP}_\varnothing=1$. We think of $\mathrm{OP}_J$ as the formal expansion of the product of all terms $(1+X_j)\in 1+k\langle\langle\underline{J}\rangle\rangle_0$, each appearing in definite order according to the ordering on $J$. We also set $\mathrm{op}_J:=\mathrm{OP}_J-1\in k\langle\langle\underline{J}\rangle\rangle_0$.

**Lemma 4.1.** *Let $m\in\mathbb{N}$ and $J=(\{0,\ldots,m-1\},<)$ where $0<\cdots<m-1$. Then*

$$\mathrm{OP}_J=(1+X_0)\cdots(1+X_{m-1}).$$

**Proof.** Immediate by definition of $\mathrm{OP}_J$. □

**Proposition 4.2.** *Let $I$ be a linearly ordered set and assume that $J=\coprod_I(J_i)_{i\in\underline{I}}$. Then $(\mathrm{op}_{J_i})_{i\in\underline{I}}$ is summable in $k\langle\langle\underline{J}\rangle\rangle_0$ and we have*

$$\mathrm{ev}_{(\mathrm{op}_{J_i})_{i\in\underline{I}}}(\mathrm{OP}_I)=\mathrm{OP}_J.$$

**Proof.** That $(\mathrm{op}_{J_i})_{i\in\underline{I}}$ is summable follows from the fact that the sets $J_i, i\in I$ are pairwise disjoint. Write $P:=\mathrm{ev}_{(\mathrm{op}_{J_i})_{i\in\underline{I}}}(\mathrm{OP}_{J_i})$. For $i\in\underline{I}$, we have $\operatorname{supp}\mathrm{op}_{j_i}\subseteq J_i^\star$. By (4.1), the set $\operatorname{supp}P$ is contained in

$$\bigcup_{(i_1,\ldots,i_n)\in I^\star}(\operatorname{supp}\mathrm{op}_{J_{i_1}}):\cdots:(\operatorname{supp}\mathrm{op}_{J_{i_n}})\subseteq\bigcup_{(i_1,\ldots,i_n)\in I^\star}J_{i_1}^\star:\cdots:J_{i_n}^\star\subseteq J^\star.$$

So $\operatorname{supp}P\subseteq J^\star=\operatorname{supp}\mathrm{OP}_J$.

Therefore it suffices in order to conclude to show that for $\theta \in J^\star$, we have $P(\theta)=1$. Let $\theta \in J^\star$ and let $p \in \mathbb{N}$ with $\theta \in (\underline{J})^p$. If $p=0$, then $\theta=\varnothing$ and $P(\varnothing)=1$ by definition. Assume now that $p>0$. We have

$$\begin{aligned} P(\theta) &= \sum_{(i_1,\ldots,i_n)\in I^\star} \mathrm{OP}_I(i_1,\ldots,i_n)\,(\mathrm{op}_{J_{i_1}}\cdots \mathrm{op}_{J_{i_n}})(\theta) \\ &= \sum_{(i_1,\ldots,i_n)\in I^\star} (\mathrm{op}_{J_{i_1}}\cdots \mathrm{op}_{J_{i_n}})(\theta). \end{aligned}$$

For $(i_1,\ldots,i_n)\in I^\star$, if the product

$$(\mathrm{op}_{J_{i_1}}\cdots \mathrm{op}_{J_{i_n}})(\theta) = \sum_{\theta_1:\cdots:\theta_n=\theta} \mathrm{op}_{J_{i_1}}(\theta_1)\cdots \mathrm{op}_{J_{i_n}}(\theta_n)$$

is non-zero, then one of its summands must be non-zero, so there must exist $\theta_1 \in J_{i_1}^\star \setminus \{\varnothing\}$,..., $\theta_n \in J_{i_n}^\star \setminus \{\varnothing\}$ with $\theta=\theta_1:\cdots:\theta_n$, in which case $\mathrm{op}_{J_{i_1}}(\theta_1)\cdots \mathrm{op}_{J_{i_n}}(\theta_n)=1\cdots 1=1$. In view of (4.1), there is a unique such $(\theta_1,\ldots,\theta_n)$. Thus the sum $\sum_{(i_1,\ldots,i_n)\in I^\star}(\mathrm{op}_{J_{i_1}}\cdots \mathrm{op}_{J_{i_n}})(\theta)$ has exactly one non-zero term, which is one, i.e. $P(\theta)=1$. We conclude that $P=\mathrm{OP}_J$. □

**Lemma 4.3.** *Let $\varsigma \notin \bigcup \underline{J}$, let $\underline{J}\times\{\varsigma\}$ be an isomorphic copy of $J$ under $j\mapsto (j,\varsigma)$ and write $\underline{J'}=\underline{J}\sqcup(\underline{J}\times\{\varsigma\})$. For all $j\in\underline{J}$, write*

$$C_J(j) := \mathrm{OP}_{j^+}^{-1}\cdot(1+X_{(j,\varsigma)})\cdot \mathrm{OP}_{j^+} \in 1+k\langle\!\langle \underline{J'}\rangle\!\rangle_0.$$

*Then $(C_J(j)-1)_{i\in\underline{j}}$ is summable in $k\langle\!\langle \underline{J'}\rangle\!\rangle_0$, and*

$$\mathrm{ev}_{((1+X_j)\cdot(1+X_{(j,\varsigma)})-1)_{j\in\underline{J}}}(\mathrm{OP}_J) = \mathrm{OP}_J\cdot \mathrm{ev}_{(C_J(j)-1)_{i\in\underline{j}}}(\mathrm{OP}_J).$$

**Proof.** For $j\in\underline{J}$, we have $C_J(j)-1=\mathrm{OP}_{j^+}^{-1}\cdot X_{(j,\varsigma)}\cdot \mathrm{OP}_{j^+}$. Since $\mathrm{OP}_{j^+}\in k\langle\!\langle \underline{J}\rangle\!\rangle$ and $(j,\varsigma)\notin\underline{J}$, any two $C_J(j)-1$ and $C_J(j')-1$ for distinct $j,j'\in\underline{J}$ have disjoint supports. It follows that $(C_J(j)-1)_{i\in\underline{j}}$ is summable. As in Proposition 4.2, it suffices to prove the result for finite $\underline{J}$. In that case, we have $C_J(j)=f[g]$ for $f=(1+X_j)_{j\in\underline{J}}$ and $g=(1+X_{(j,\varsigma)})_{j\in\underline{J}}$, so this follows from Lemma 4.1. □

**Lemma 4.4.** *Let $J,I$ be sets, let $\varepsilon\in k\langle\!\langle J\rangle\!\rangle_0$ and let $Q: I\longrightarrow k\langle\!\langle J\rangle\!\rangle_0$ be a summable family and write $(1+\varepsilon)\cdot Q\cdot(1+\varepsilon)^{-1}$ for the family $((1+\varepsilon)\cdot Q(i)\cdot(1+\varepsilon)^{-1})_{i\in I}$. Then for all $P\in k\langle\!\langle J\rangle\!\rangle$, we have*

$$\mathrm{ev}_{(1+\varepsilon)\cdot Q\cdot(1+\varepsilon)^{-1}}(P) = (1+\varepsilon)\cdot \mathrm{ev}_Q(P)\cdot(1+\varepsilon)^{-1}.$$

**Proof.** By Proposition 3.2, the conjugacy function

$$\begin{aligned} k\langle\!\langle J\rangle\!\rangle &\longrightarrow k\langle\!\langle J\rangle\!\rangle \\ P &\longmapsto (1+\varepsilon)\cdot P\cdot(1+\varepsilon)^{-1} \end{aligned}$$

is a strongly linear automorphism of $k\langle\!\langle J\rangle\!\rangle$ that sends $k\langle\!\langle J\rangle\!\rangle_0$ into $k\langle\!\langle J\rangle\!\rangle_0$. The result follows. □

**Corollary 4.5.** *Let $J$ be a linearly ordered set, let $\varsigma\notin\underline{J}$ and let $J':=\underline{J}\sqcup\{\varsigma\}$. Write $P$ for the family $\underline{J}\longrightarrow k\langle\!\langle J'\rangle\!\rangle_0$ with $P(j):=(1+X_\varsigma)\cdot X_j\cdot(1+X_\varsigma)^{-1}$ for all $j\in\underline{J}$. Then*

$$\mathrm{ev}_Q(\mathrm{OP}_J) = (1+X_\varsigma)\cdot \mathrm{OP}_J\cdot(1+X_\varsigma)^{-1}.$$

**Proof.** Apply Lemma 4.4 to $\varepsilon=X_\varsigma$ and $Q(j)=X_j$ for all $j\in\underline{J}$. □

**Lemma 4.6.** *Let* $I=(\underline{I},<)$ *be a linearly ordered set. We define a family* $R:\underline{I}\longrightarrow k\langle\langle\underline{I}\rangle\rangle$ *by*

$$\forall i\in\underline{I}, R(i):=\sum_{n\in\mathbb{N}}(-1)^n X_i^n=(1+X_i)^{-1}.$$

*We have*

$$\mathrm{ev}_{R-1}(\mathrm{OP}_{I^*})=\mathrm{OP}_I^{-1}.$$

**Proof.** Let $n\in\mathbb{N}$ and $\theta=(j_1,\ldots,j_n)\in(\underline{I})^n$. Let us show that $\mathrm{ev}_{R-1}(\mathrm{OP}_{I^*})(\theta)=\mathrm{OP}_I^{-1}(\theta)$. Write $S:=(\{j_1,\ldots,j_p\},<)$ and $R'=R\upharpoonright\{j_1,\ldots,j_n\}$. We have

$$\begin{aligned}\mathrm{ev}_{R-1}(\mathrm{OP}_{I^*})(\theta) &= \mathrm{ev}_{R'-1}(\mathrm{OP}_{S^*})(\theta)\quad\text{and}\\ \mathrm{OP}_I^{-1}(\theta) &= \mathrm{OP}_S^{-1}(\theta),\end{aligned}$$

so we may assume that $\underline{I}$ is finite (hence also $R'=R$). Write $I=\{i_1,\ldots,i_m\}$ where $m\leqslant n$ and $i_1<\cdots<i_n$. By Lemma 4.1, we have $\mathrm{OP}_{I^*}=(1+X_{i_n})\cdots(1+X_{i_1})$ and $\mathrm{OP}_I=(1+X_{i_1})\cdots(1+X_{i_n})$, so

$$\mathrm{ev}_{R-1}(\mathrm{OP}_{I^*})=R(i_n)\cdots R(i_1)=(1+X_{i_n})^{-1}\cdots(1+X_{i_1})^{-1}=(\mathrm{OP}_I)^{-1}.$$

In particular $\mathrm{ev}_{R-1}(\mathrm{OP}_{I^*})(\theta)=\mathrm{OP}_I^{-1}(\theta)$. This concludes the proof. □

**Lemma 4.7.** *Suppose that* $k$ *has characteristic zero. Let* $I=(\underline{I},<)$ *be a linearly ordered set, and let* $L$ *be a summability Lie subalgebra of* $k\langle\langle\underline{I}\rangle\rangle_0$ *containing the set* $\{\log(1+X_i):i\in\underline{I}\}$*. Then* $\log(\mathrm{OP}_I)\in L$*.*

**Proof.** We prove this by induction on the cardinal $\kappa:=|\underline{I}|$. In view of Lemma 4.1 and by [4, Proposition 2.11], the result holds if $\kappa\leqslant 2$. By induction, it holds if $\kappa<\omega$. Let $\kappa$ such that the result holds for all cardinals $\lambda<\kappa$. Let $f:\kappa\longrightarrow\underline{I}$ be a bijection, and write $J_\alpha:=\{f(\beta):\beta<\alpha\}$ for each $\alpha\leqslant\kappa$. So $J_\alpha$ has cardinality $<\kappa$ whenever $\alpha<\kappa$, and $J_\kappa=\underline{I}$.

Let $J\subseteq\underline{I}$. Note that $\rho_J(\mathrm{op}_I)=(\sum_{\theta\in I^\star\cap(J,<)^\star}X_\theta)-1=\mathrm{op}_{(J,<)}$. Thus

$$\begin{aligned}\rho_J(\log(\mathrm{OP}_I)) &= \rho_J\left(\sum_{m>0}\frac{(-1)^{m+1}}{m}\mathrm{op}_I^m\right)\\ &= \sum_{m>0}\frac{(-1)^{m+1}}{m}(\rho_J(\mathrm{op}_I))^m\\ &= \sum_{m>0}\frac{(-1)^{m+1}}{m}\mathrm{op}_{(J,<)}^m\\ &= \log(\mathrm{OP}_{(J,<)}).\end{aligned}$$

Therefore $\log(\mathrm{OP}_I)\in\log(\mathrm{OP}_{(J,<)})+\mathrm{Ker}(\rho_J)$ for each $J\subseteq\underline{I}$. So for each $\alpha<\kappa$, we have

$$\log(\mathrm{OP}_I)-\log(\mathrm{OP}_{(J_\alpha,<)})\in\mathrm{Ker}(\rho_{J_\alpha}).\tag{4.2}$$

Write $Q_\alpha:=\log(\mathrm{OP}_{(J_{\alpha+1},<)})-\log(\mathrm{OP}_{(J_\alpha,<)})$ for each $\alpha<\kappa$. Our induction hypothesis for $L$ gives $\log(\mathrm{OP}_{(J_{\alpha+1},<)}),\log(\mathrm{OP}_{(J_\alpha,<)})\in L$, whence $Q_\alpha\in L$. Furthermore [4, Proposition 2.11] implies that $Q_\alpha$ lies in the summability subspace of $k\langle\langle J_{\alpha+1}\rangle\rangle$ generated by $\log(1+X_{f(\alpha)})$, Lie brackets of $\log(1+X_{f(\alpha)})$ and elements in $k\langle\langle J_\alpha\rangle\rangle$. It follows since $\mathrm{supp}\log(1+X_{f(\alpha)})\subseteq J_{\alpha+1}^\star\setminus J_\alpha^\star$ and in view of (3.1) that

$$\mathrm{supp}\,Q_\alpha\subseteq J_{\alpha+1}^\star\setminus J_\alpha^\star.\tag{4.3}$$

In particular, the family $(Q_\alpha)_{\alpha<\kappa}$ is summable in $k\langle\langle \underline{I}\rangle\rangle$. We deduce since $L$ is a summability subspace that $\sum_{\alpha<\kappa} Q_\alpha \in L$. We claim that $\sum_{\alpha<\kappa} Q_\alpha = \log(\mathrm{OP}_I)$, and we will prove by induction on $\alpha \leqslant \kappa$ that $\sum_{\beta<\alpha} Q_\beta = \log(\mathrm{OP}_{(J_\alpha,<)})$. This holds for $\alpha = 0$ since $\log(1) = 0$. Let $\alpha \leqslant \kappa$ such that the result holds for all $\beta < \alpha$. If $\alpha = \beta + 1$ is a successor, then we immediately get the desired identity. Assume that $\alpha$ is a non-zero limit.

For $\theta \in J_\alpha^\star$, there is a least $\gamma \leqslant \alpha$ with $\theta \in J_\gamma^\star$ (since $\theta$ is a finite word). We have $\theta \in J_\gamma^\star \setminus J_\eta^\star$ for all $\eta < \gamma$ by minimality, and $\theta \notin \mathrm{Ker}(\rho_{J_\gamma})$. So

$$\begin{aligned} \log(\mathrm{OP}_{(J_\alpha,<)})(\theta) &= \log(\mathrm{OP}_{(J_\gamma,<)})(\theta) && \text{(by (4.2))} \\ &= \sum_{\eta<\gamma} Q_\gamma(\theta) \\ &= \sum_{\beta<\alpha} Q_\beta(\theta). && \text{(by (4.3))} \end{aligned}$$

We deduce by induction that $\log(\mathrm{OP}_I) = \log(\mathrm{OP}_{(J_\kappa,<)}) = \sum_{\beta<\kappa} Q_\beta \in L$. □

## 4.2 Ordered products in $1 + \mathfrak{m}$

In this section, we fix a summability algebra $A = k + \mathfrak{m}$ with evaluations. We also fix a linearly ordered set $I$ and an element $\varsigma \notin \underline{I} \cup (\bigcup \underline{I})$.

A family $f : \underline{I} \longrightarrow 1 + \mathfrak{m}$ is said multipliable if $f - 1 : \underline{I} \longrightarrow \mathfrak{m}$ is summable. We then define the *ordered product* $\Pi_I f$ of $f$ as

$$\Pi_I f := \mathrm{ev}_{f-1}(\mathrm{OP}_I) \in 1 + \mathfrak{m}.$$

We fix a multipliable family $f : \underline{I} \longrightarrow 1 + \mathfrak{m}$.

**Proposition 4.8.** *The set* $\mathrm{dom}\,\Pi_I$ *is a subgroup of* $(1 + \mathfrak{m})^{\underline{I}}$.

**Proof.** Let $g : \underline{I} \longrightarrow 1 + \mathfrak{m}$ be multipliable. By **SA** in $(A, \Sigma)$, the family $(f - 1) \cdot (g - 1)$ is summable, so $f \cdot g - 1 = (f - 1) + (g - 1) + (f - 1) \cdot (g - 1)$ is summable. So $f \cdot g \in \mathrm{dom}\,\Pi_I$. Since $(A, \Sigma)$ has evaluations, the family $((-1)^k\,(f(i) - 1)^k)_{i \in \underline{I} \wedge k > 0}$ is summable. By **SS3**, so is $(\sum_{k>0} (-1)^k\,(f(i) - 1)^k)_{i \in \underline{I}} = (f(i)^{-1} - 1)_{i \in \underline{I}}$. Therefore $f^{-1} \in \mathrm{dom}\,\Pi_I$. So $\mathrm{dom}\,\Pi_I$ is a subgroup of $(1 + \mathfrak{m})^{\underline{I}}$. □

**Proposition 4.9.** *Let* $g : \underline{I} \longrightarrow 1 + \mathfrak{m}$ *be multipliable. Consider the family*

$$f[g] := (\Pi_{(i^+)^*}(f^{-1} \upharpoonright (\underline{i^+})) \cdot g(i) \cdot \Pi_{i^+}(f \upharpoonright (\underline{i^+})))_{i \in \underline{I}}.$$

*Then* $f[g]$ *is multipliable, with* $\Pi_I(f \cdot g) = (\Pi_I f) \cdot (\Pi_I f[g])$.

**Proof.** In the notations of Lemma 4.3, let $\varphi : \underline{I} \longrightarrow \underline{I} \times \{\varsigma\}; i \mapsto (i, \varsigma)$ be the natural isomorphism and write $g' = g \circ \varphi$. So $g' - 1$ is summable. We have

$$\mathrm{ev}_{((1+X_i)\cdot(1+X_{(i,\varsigma)})-1)_{i\in\underline{I}}}(\mathrm{OP}_I) = \mathrm{OP}_I \cdot \mathrm{ev}_{(C_I(i)-1)_{i\in\underline{I}}}(\mathrm{OP}_I).$$

Therefore

$$\begin{aligned} \Pi_I(f \cdot g) &= \mathrm{ev}_{(f(i)\cdot g(i))_{i\in\underline{I}}}(\mathrm{OP}_I) \\ &= \mathrm{ev}_{\mathrm{ev}_{f\sqcup g'}(((1+X_i)\cdot(1+X_{(i,\varsigma)})-1))_{i\in\underline{I}}}(\mathrm{OP}_I) \\ &= \mathrm{ev}_{f\sqcup g'}(\mathrm{ev}_{((1+X_i)\cdot(1+X_{(i,\varsigma)})-1))_{i\in\underline{I}}})(\mathrm{OP}_I) && \text{(by Proposition 3.4)} \\ &= \mathrm{ev}_{f\sqcup g'}(\mathrm{OP}_I \cdot \mathrm{ev}_{(C_I(i)-1)_{i\in\underline{i}}}(\mathrm{OP}_I)) \\ &= \mathrm{ev}_{f\sqcup g'}(\mathrm{OP}_I) \cdot \mathrm{ev}_{f\sqcup g'}(\mathrm{ev}_{(C_I(i)-1)_{i\in\underline{I}}}(\mathrm{OP}_I))) \\ &= (\Pi_I f) \cdot \mathrm{ev}_{(\mathrm{ev}_{f\sqcup g'}(C_I(i)-1))_{i\in\underline{I}}}(\mathrm{OP}_I) && \text{(by Proposition 3.4)} \\ &= (\Pi_I f) \cdot \mathrm{ev}_{((f[g](i)-1))_{i\in\underline{I}}}(\mathrm{OP}_I) \\ &= (\Pi_I f) \cdot (\Pi_I f[g]). \end{aligned}$$

This concludes the proof. $\square$

**Proposition 4.10.** *Assume that $f$ has finite support* $\operatorname{supp} f = \{i_1, \ldots, i_n\}$ *where* $i_1 < \cdots < i_n$. *Then* $\Pi_I f = f(i_1) \cdots f(i_n)$.

**Proof.** By Lemma 4.1, we have $\Pi_I f = (1 + f(i_1) - 1) \cdots (1 + f(i_n) - 1) = f(i_1) \cdots f(i_n)$. $\square$

**Proposition 4.11.** *Let $J$ be a linearly ordered set and let $\varphi: J \longrightarrow I$ be an isomorphism. Then $f \circ \varphi$ is multipliable, and $\Pi_J (f \circ \varphi) = \Pi_I f$.*

**Proof.** The function $\varphi$ is bijective, so **SS2** in $(A, \Sigma)$ implies that $f \circ \varphi - 1$ is summable. The map $\varphi$ induces a bijection $\varphi^\star: J^\star \longrightarrow I^\star; (j_1, \ldots, j_n) \longmapsto (\varphi(j_1), \ldots, \varphi(j_n))$. Write

$$\begin{aligned} F &:= ((f(i_1) - 1) \cdots (f(i_n) - 1))_{(i_1, \ldots, i_n) \in I^\star} \quad \text{and} \\ G &:= ((f \circ \varphi(j_1) - 1) \cdots (f \circ \varphi(j_n) - 1))_{(j_1, \ldots, j_n) \in J^\star}. \end{aligned}$$

We have $G = F \circ \varphi^\star$ where $\Pi_I f = \Sigma_{I^\star} F$ and $\Pi_J (f \circ \varphi) = \Sigma_{J^\star} G$, so **SS2** in $(A, \Sigma)$ gives $\Pi_I f = \Pi_J (f \circ \varphi)$. $\square$

**Proposition 4.12.** *Let $J$ be a linearly ordered set, and let $(I_j)_{j \in J}$ be a family of linearly ordered subsets with $J = \coprod_{j \in J} I_j$. Let $g: \underline{I} \longrightarrow 1 + \mathfrak{m}$ be multipliable. For each $j \in \underline{J}$, write $g_j := g \upharpoonright \underline{I_j}$. Then*

*a)* $g_j$ *is multipliable,*

*b)* $(\Pi_{I_j} g_j)_{j \in \underline{J}} - 1$ *is multipliable, and*

*c)* $\Pi_J (\Pi_{I_j} g_j)_{j \in J} = \Pi_I g$.

**Proof.** That $g_j - 1$ is summable follows from **SS3a** in $(A, \Sigma)$. Since $(A, \Sigma)$ is has evaluations, the family $(\operatorname{op}_I(i_1, \ldots, i_n)\, (g(i_1) - 1) \cdots (g(i_n) - 1))_{(i_1, \ldots, i_n) \in I^\star}$ is summable. Let $I' := \bigsqcup_{j \in \underline{J}} I_j^\star$, so $I' \subseteq I^\star$. The family $(\operatorname{op}_I(i_1, \ldots, i_n)\, (g(i_1) - 1) \cdots (g(i_n) - 1))_{(i_1, \ldots, i_n) \in I'}$ is also summable by **SS3a**. By **SS3b** in $(A, \Sigma)$ for the disjoint union $I' = \bigsqcup_{j \in \underline{J}} I_j^\star$, the family

$$\left( \sum_{(i_1, \ldots, i_n) \in I_j} \operatorname{op}_I(i_1, \ldots, i_n)\, (g(i_1) - 1) \cdots (g(i_n) - 1) \right)_{j \in \underline{J}},$$

is summable. But this is none other than $(\Pi_I g_j)_{j \in \underline{J}}$. So *b*) holds. Finally, we have

$$\begin{aligned} \Pi_I (\Pi_{I_j} g_j)_{j \in \underline{J}} &= \operatorname{ev}_{((\Pi_{I_j} g_j) - 1)_{j \in \underline{J}}}(\mathrm{OP}_J) \\ &= \operatorname{ev}_{(\operatorname{ev}_{g_j - 1}(\operatorname{op}_{I_j}))_{j \in J}}(\mathrm{OP}_J) \\ &= \operatorname{ev}_{(\operatorname{ev}_{g - 1}(\operatorname{op}_{I_j}))_{j \in J}}(\mathrm{OP}_J) && \text{(since each } \operatorname{op}_{I_j} \in k\langle\langle \underline{I_j} \rangle\rangle\text{)} \\ &= \operatorname{ev}_{g - 1}(\operatorname{ev}_{(\operatorname{op}_{I_j})_{j \in J}}(\mathrm{OP}_J)) && \text{(by Proposition 3.4)} \\ &= \operatorname{ev}_{g - 1}(\mathrm{OP}_I) && \text{(by Proposition 4.2)} \\ &= \Pi_I g. \end{aligned}$$

This proves *c*). $\square$

**Proposition 4.13.** *Let $\varepsilon \in \mathfrak{m}$. Then the family*

$$(1 + \varepsilon) \cdot f \cdot (1 + \varepsilon)^{-1} = ((1 + \varepsilon) \cdot f(i) \cdot (1 + \varepsilon)^{-1})_{i \in \underline{I}}$$

*is multipliable and* $\Pi_I ((1 + \varepsilon) \cdot f \cdot (1 + \varepsilon)^{-1}) = (1 + \varepsilon) \cdot (\Pi_I f) \cdot (1 + \varepsilon)^{-1}$.

**Proof.** That $(1+\varepsilon)\cdot f\cdot(1+\varepsilon)^{-1}$ is multipliable follows from the fact that $(A,\Sigma)$ is a strong algebra. That it takes values in $\mathfrak{m}$ follows from the fact that $\mathfrak{m}$ is a bilateral ideal in $A$. Consider $I' := \underline{I} \sqcup \{\varsigma\}$ as in Corollary 4.5 and write $f'$ for the family $f-1 \sqcup \{(\varsigma, \varepsilon)\}: I' \longrightarrow 1+\mathfrak{m}$. We have

$$\begin{aligned}
\Pi_I((1+\varepsilon)\cdot f\cdot(1+\varepsilon)^{-1}) &= \mathrm{ev}_{(1+\varepsilon)\cdot f\cdot(1+\varepsilon)^{-1}-1}(\mathrm{OP}_I)\\
&= \mathrm{ev}_{(\mathrm{ev}_{f'}((1+X_\varsigma)\cdot X_i\cdot(1+X_\varsigma)^{-1}))_{i\in I}}(\mathrm{OP}_I)\\
&= \mathrm{ev}_{f'}(\mathrm{ev}_{(1+X_\varsigma)\cdot X_i\cdot(1+X_\varsigma)^{-1}}(\mathrm{OP}_I)) && \text{(by Corollary 4.5)}\\
&= \mathrm{ev}_{f'}((1+X_\varsigma)\cdot \mathrm{OP}_I\cdot(1+X_\varsigma)^{-1}) && \text{(by Proposition 3.4)}\\
&= \mathrm{ev}_{f'}(1+X_\varsigma)\cdot \mathrm{ev}_{f'}(\mathrm{OP}_I)\cdot \mathrm{ev}_{f'}(1+X_\varsigma)^{-1}\\
&= (1+\varepsilon)\cdot(\Pi_I f)\cdot(1+\varepsilon)^{-1}.
\end{aligned}$$

This concludes the proof. □

**Proposition 4.14.** *The family $f^{-1} = (f(i)^{-1})_{i\in\underline{I}}$ is multipliable and $\Pi_{I^*}(f^{-1}) = (\Pi_I f)^{-1}$.*

**Proof.** We define a family $P: I \longrightarrow k\langle\langle\underline{I}\rangle\rangle$ by setting $R(i) := (1+X_i)^{-1}$ for each $i\in\underline{I}$. Set $Q := \sum_{i\in\underline{I}} R(i) \in k\langle\langle\underline{I}\rangle\rangle$. The family $(Q(i_1,\ldots,i_n)\,(f(i_1)-1)\cdots(f(i_n)-1))_{(i_1,\ldots,i_n)\in\underline{I}^\star}$ is summable. Now by **SS3** in $(A,\Sigma)$, the family $f^{-1}-1 = (\mathrm{ev}_{f-1}(R(i)))_{i\in\underline{I}}$, i.e.

$$\left(\sum_{(i_1,\ldots,i_n)\in\{i\}^\star} Q(i_1,\ldots,i_n)\,(f(i_1)-1)\cdots(f(i_n)-1)\right)_{i\in\underline{I}}$$

is summable. We have

$$\begin{aligned}
\Pi_{I^*} f^{-1} &= \mathrm{ev}_{f^{-1}-1}(\mathrm{OP}_{I^*})\\
&= \mathrm{ev}_{(\mathrm{ev}_{f-1}(R(i))-1)_{i\in\underline{I}}}(\mathrm{OP}_{I^*})\\
&= \mathrm{ev}_{f-1}(\mathrm{ev}_{R-1}(\mathrm{OP}_{I^*})) && \text{(by Proposition 3.4)}\\
&= \mathrm{ev}_{f-1}(\mathrm{OP}_I^{-1}) && \text{(by Lemma 4.6)}\\
&= \mathrm{ev}_{f-1}(\mathrm{OP}_I)^{-1}\\
&= (\Pi_I f)^{-1}.
\end{aligned}$$

This concludes the proof. □

## 4.3 The multipliability group $1+\mathfrak{m}$

**Theorem 4.15.** *Let $(A,\Sigma)$ be a summability algebra with evaluations, and write $\mathfrak{m}$ for its maximal ideal. For each linearly ordered set I, let $\Pi_I$ be the ordered product partial function on $(1+\mathfrak{m})^{\underline{I}}$. Then $(1+\mathfrak{m}, 1, \cdot, \Pi)$ is a full multipliability group.*

**Proof.** The validity of the axioms **MG1–MG7** follows from Propositions 4.8–4.14. The axiom **Fullness** is satisfied by definition of multipliability in $1+\mathfrak{m}$. □

**Proposition 4.16.** *Assume that $k$ has characteristic zero. Let $G\subseteq 1+\mathfrak{m}$ be a subset such that $\log(G)$ is a summability Lie subalgebra of $\mathfrak{m}$. Then $G$ is a multipliability subgroup of $1+\mathfrak{m}$.*

**Proof.** Let $a,b\in G$. By [4, Proposition 2.11], the element $\log(a)*\log(b^{-1})$ lies in all summability Lie subalgebras of $\mathfrak{m}$ containing $\log(a)$ and $\log(b^{-1}) = -\log(b)$. So $\log(a\cdot b^{-1}) = \log(a)*\log(b^{-1})\in\log(G)$, whence $a\cdot b^{-1}\in G$. Therefore $G$ is a subgroup of $1+\mathfrak{m}$.

Let $I = (\underline{I}, <)$ be a linearly ordered set and let $(g_i)_{i \in \underline{I}}$ be a family in $G$ which is multipliable in $1+\mathfrak{m}$. We want to show that its product $g = \prod_I (g_i)_{i \in \underline{I}}$ lies in $G$. Since log is injective, it suffices to show that $\log(g) \in \log(G)$. Let $\delta = (g_i - 1)_{i \in \underline{I}}$, so $\delta$ is summable in $\mathfrak{m}$. We have

$$\begin{aligned} \log(g) &= \log(\mathrm{ev}_\delta(\mathrm{OP}_I)) && \\ &= \mathrm{ev}_{\mathrm{ev}_\delta(\mathrm{OP}_I)-1}(\log(1+X_\bullet)) && \text{(by definition of log)} \\ &= \mathrm{ev}_\delta(\mathrm{ev}_{\mathrm{OP}_I - 1}(\log(1+X_\bullet))) && \text{(by Proposition 3.4)} \\ &= \mathrm{ev}_\delta(\log(\mathrm{OP}_I)). && \text{(by definition of log)} \end{aligned}$$

By Lemma 4.7, the formal series $\log(\mathrm{OP}_I)$ lies in all summability Lie subalgebras of $k\langle\langle \underline{I} \rangle\rangle_0$ containing $\{\log(1+X_i) : i \in \underline{I}\}$. Consider the subspace $L := \mathrm{ev}_\delta^{-1}(\log(G))$ of $k\langle\langle \underline{I} \rangle\rangle_0$. Since $\mathrm{ev}_\delta : k\langle\langle \underline{I} \rangle\rangle_0 \longrightarrow \mathfrak{m}$ is a morphism of algebras, this is a Lie subalgebra of $k\langle\langle \underline{I} \rangle\rangle_0$. It contains $\{\log(1+X_i) : i \in \underline{I}\}$ since $\mathrm{ev}_\delta(\log(1+X_i)) = \log(g_i) \in \log(G)$ for each $i \in \underline{I}$. Lastly, if $J$ is a set and $(a_j)_{j \in J} \in L^J$ is summable in $k\langle\langle \underline{I} \rangle\rangle_0$, then we have $\mathrm{ev}_\delta(\sum_{j \in J} a_j) = \sum_{j \in J} \mathrm{ev}_\delta(a_j) \in \log(G)$. This shows that $L$ is a summability subspace of $k\langle\langle \underline{I} \rangle\rangle_0$. So $\log(\mathrm{OP}_I) \in L$, whence $\log(g) \in \log(G)$. □

**Remark 4.17.** Suppose that $k$ has characteristic 0. Then the axiom **UF** translates into a similar axiom for $1+\mathfrak{m}$: if $(f_i)_{i \in I}$ is multipliable and $(c_i)_{i \in I}$ is a family of functions $c_i : X_i \longrightarrow k$ with finite domains $X_i$, then $\left(f_i^{[c_i(x)]}\right)_{i \in I \wedge x \in X_i}$ is also multipliable. We believe that full multipliability $k$-groups satisfying this and in which sufficiently many group-like elements of $k\langle\langle I \rangle\rangle$ (see [24]) can be evaluated at multipliable functions are in correspondence with summability $k$-algebras with evaluations, *via* the Lie correspondence of [4]. The equivalence should be an equivalence of suitable categories, based on the Mal'cev-Lazard-Stewart correspondence [18, 17, 25]. An established case of this correspondence will be implicitely at play in our applications to flat transseries in Section 6.

# 5 Chain rules in algebras of Noetherian series

## 5.1 A formal setting for chain rules

Let $I = (\underline{I}, <)$ be a non-empty linearly ordered set. Let $\varsigma \notin \underline{I} \cup \bigcup \underline{I}$, and consider the set

$$I_\varsigma := \underline{I} \sqcup \{\varsigma\} \sqcup (\underline{I} \times \{\varsigma\}).$$

We have a summability algebra $k\langle\langle I_\varsigma \rangle\rangle$. Let $\mathcal{R}(\underline{I})$ be the subset of $k\langle\langle J \rangle\rangle_0$ containing

$$(1+X_i) \cdot X_\varsigma - (1+X_{(i,\varsigma)}) \cdot X_\varsigma \cdot (1+X_i), \tag{5.1}$$

for each $i \in \underline{I}$. We write $\mathfrak{q}(\underline{I})$ for the smallest summability ideal of $k\langle\langle J \rangle\rangle_0$ containing $\mathcal{R}(\underline{I})$.

Consider the family $Q_I : \underline{I} \longrightarrow 1 + k\langle\langle J \rangle\rangle_0$ given by

$$Q_I(i) := \mathrm{OP}_{i^-} \cdot (1+X_{(i,\varsigma)}) \cdot \mathrm{OP}_{i^-}^{-1}$$

where $i^- = \{j \in \underline{I} : j < i\}$ for each $i \in \underline{I}$. The occurrence of $X_{(i,\varsigma)}$ in the definition of $Q_I(i)$ implies that $\mathrm{supp}\,(Q_I(i) - 1) \cap \mathrm{supp}\,(Q_I(j) - 1) = \varnothing$ whenever $i, j \in \underline{I}$ are distinct. Therefore $Q_I - 1 : \underline{I} \longrightarrow k\langle\langle J \rangle\rangle_0$ is summable. As a consequence, the following ordered product is well-defined:

$$M_I := \mathrm{ev}_{Q_I - 1}(\mathrm{OP}_I) = \prod_I (Q_I(i))_{i \in \underline{I}} \in 1 + k\langle\langle J \rangle\rangle_0.$$

**Lemma 5.1.** *We have* $M_I \cdot X_\varsigma \cdot \mathrm{OP}_I - \mathrm{OP}_I \cdot X_\varsigma \in \mathfrak{q}(\underline{I})$.

**Proof.** We proceed in a similar way as in the proof of Lemma 4.7. Let $\kappa$ be the cardinality of $\underline{I}$ and let $f\colon \kappa \longrightarrow \underline{I}$ be a bijection. For all $\alpha \leqslant \kappa$, let $I_\alpha$ denote the subset $I_\alpha := \{f(\beta): \beta < \alpha\} \subseteq \underline{I}$. Write

$$\begin{aligned} R_\alpha &:= M_{(I_\alpha,<)} \cdot X_\varsigma \cdot \mathrm{OP}_{(I_\alpha,<)} - \mathrm{OP}_{(I_\alpha,<)} \cdot X_\varsigma \qquad \text{and} \\ D_{\alpha,\beta} &:= R_\beta - R_\alpha \end{aligned}$$

for each $\alpha, \beta \leqslant \kappa$. Note that $\operatorname{supp} R_\alpha \subseteq ((I_\alpha)_\varsigma)^\star$ and $\operatorname{supp} D_{\alpha,\beta} \subseteq ((I_{\max(\alpha,\beta)})_\varsigma)^\star$.

For each $\alpha \leqslant \beta \leqslant \kappa$, we can write $\mathrm{OP}_{(I_\beta,<)} = \mathrm{OP}_{(I_\alpha,<)} + N_{\alpha,\beta}$ and $M_{(I_\beta,<)} = M_{(I_\alpha,<)} + T_{\alpha,\beta}$ where $\operatorname{supp} N_\alpha, \operatorname{supp} T_\alpha \subseteq ((I_\beta)_\varsigma)^\star \setminus ((I_\alpha)_\varsigma)^\star$. Therefore the series $D_{\beta,\alpha} = R_\beta - R_\alpha = T_{\alpha,\beta} \cdot X_\varsigma \cdot (\mathrm{OP}_{(I_\alpha,<)} + N_{\alpha,\beta}) + M_{(I_\alpha,<)} \cdot X_\varsigma \cdot N_{\alpha,\beta} - N_{\alpha,\beta} \cdot X_\varsigma$ satisfies

$$\operatorname{supp} D_{\alpha,\beta} \subseteq ((I_\beta)_\varsigma)^\star \setminus ((I_\alpha)_\varsigma)^\star. \tag{5.2}$$

In particular, the supports of elements of the family $(D_{\beta,\beta+1})_{\beta<\kappa}$ are pairwise disjoint, so that family is summable in $k\langle\langle I_\varsigma \rangle\rangle$.

We first prove by induction on $\alpha \leqslant \kappa$ that we have

$$R_\alpha = \sum_{\beta<\alpha} D_{\beta,\beta+1}. \tag{5.3}$$

Note that $I_0 = \varnothing$, so

$$R_0 = M_{(I_0,<)} \cdot X_\varsigma \cdot \mathrm{OP}_{(I_0,<)} - \mathrm{OP}_{(I_0,<)} \cdot X_\varsigma = 1 \cdot X_\varsigma \cdot 1 - 1 \cdot X_\varsigma = 0. \tag{5.4}$$

So the identity (5.3) holds at $\alpha = 0$. Let $\alpha \leqslant \kappa$ such that it holds for all $\gamma < \alpha$. If $\alpha$ is a successor, it follows immediately that it holds at $\alpha$. Assume that $\alpha$ is a non-zero limit. Consider a $\theta \in (I_\alpha)_\varsigma$. There is a $\gamma < \alpha$ such that $\theta \in (I_\gamma)_\varsigma$. We deduce in view of (5.2) and using the induction hypothesis that

$$R_\alpha(\theta) = R_\gamma(\theta) = \sum_{\beta<\gamma} D_{\beta,\beta+1}(\theta) = \sum_{\beta<\alpha} D_{\beta,\beta+1}(\theta),$$

so $R_\alpha = \sum_{\beta<\alpha} D_{\beta,\beta+1}$.

We next claim and we prove by induction on $\kappa$ that $R_\kappa \in \mathfrak{q}(\underline{I})$. Assume first that $\kappa < \omega$, so $\underline{I}$ is finite. Write $\underline{I} = \{i_1, \ldots, i_\kappa\}$ with $i_1 < \cdots < i_\kappa$. If $\kappa = 0$, then the result follows from (5.4). Assume that $\kappa > 0$ and set $\underline{I}' := \{i_1, \ldots, i_\kappa\}$ and $I' := (\underline{I}', <)$.

Recall by Lemma 4.1 that $\mathrm{OP}_I = (1 + X_{i_1}) \cdot \mathrm{OP}_{I'}$ whereas $\mathrm{OP}_{i_m^-} = (1 + X_{i_1}) \cdots (1 + X_{i_{m-1}})$ for all $m \in \{1, \ldots, n\}$. We thus see that

$$M_I = M_{I'} \cdot \mathrm{OP}_{I'} \cdot (1 + X_{(i_n,\varsigma)}) \cdot \mathrm{OP}_{I'}^{-1}. \tag{5.5}$$

Define

$$\begin{aligned} a &:= (1 + X_{(i_n,\varsigma)}) \cdot (1 + X_{i_n}) \cdot X_\varsigma - X_\varsigma \cdot (1 + X_{i_n}) \quad \text{and} \\ b &:= M_{I'} \cdot \mathrm{OP}_{I'} \cdot X_\varsigma - X_\varsigma \cdot \mathrm{OP}_{I'}, \end{aligned}$$

so $a \in \mathfrak{q}(\underline{I})$ by (5.1) and $b \in \mathfrak{q}(\underline{I}') \subseteq \mathfrak{q}(\underline{I})$ by the induction hypothesis. We have

$$\begin{aligned} X_\varsigma \cdot \mathrm{OP}_I &= X_\varsigma \cdot \mathrm{OP}_{I'} \cdot (1 + X_{i_n}) \\ &= (M_{I'} \cdot \mathrm{OP}_{I'} \cdot X_\varsigma - b) \cdot (1 + X_{i_n}) \\ &= (M_{I'} \cdot \mathrm{OP}_{I'} \cdot X_\varsigma) \cdot (1 + X_{i_n}) + d \end{aligned}$$

where $d := -b \cdot (1 + X_{i_n}) \in \mathfrak{q}(\underline{I})$. Now

$$(M_{I'} \cdot \mathrm{OP}_{I'} \cdot X_\varsigma) \cdot (1 + X_{i_n}) = M_{I'} \cdot \mathrm{OP}_{I'} \cdot (1 + X_{(i_n,\varsigma)}) \cdot (1 + X_{i_n}) \cdot X_\varsigma + d$$

where $d := M_{I'} \cdot \mathrm{OP}_{I'} \cdot a \in \mathfrak{q}(\underline{I})$. In view of (5.5), we obtain

$$\begin{aligned} M_{I'} \cdot \mathrm{OP}_{I'} \cdot (1 + X_{(i_n,\varsigma)}) \cdot (1 + X_{i_n}) \cdot X_\varsigma &= M_I \cdot \mathrm{OP}_{I'} \cdot (1 + X_{i_n}) \cdot X_\varsigma \\ &= M_I \cdot \mathrm{OP}_I \cdot X_\varsigma. \end{aligned}$$

We deduce that $R_\kappa = X_\varsigma \cdot \mathrm{OP}_I - M_I \cdot \mathrm{OP}_I \cdot X_\varsigma \in \mathfrak{q}(\underline{I})$ as claimed. By induction, the result holds whenever $\kappa$ is finite.

Now let $\kappa$ be infinite and such that the result holds for all strictly smaller cardinalities. In particular, we have $R_\gamma \in \mathfrak{q}(\underline{I_\gamma}) \subseteq \mathfrak{q}(\underline{I})$ for all $\gamma < \kappa$. It follows since $\mathfrak{q}(\underline{I})$ is a subgroup of $k\langle\langle I_\varsigma \rangle\rangle$ that $D_{\gamma,\gamma+1} \in \mathfrak{q}(\underline{I})$ for each $\gamma < \kappa$. Finally, we deduce since $\mathfrak{q}(\underline{I})$ is a summability subspace that $R_\kappa = \sum_{\gamma<\kappa} D_{\gamma,\gamma+1} \in \mathfrak{q}(\underline{I})$. This concludes the proof. □

**Lemma 5.2.** *Let $(A, \Sigma)$ be a summability algebra with evaluations, with maximal ideal $\mathfrak{m}$. Let $\varepsilon \in \mathfrak{m}$ and let $a : \underline{I} \longrightarrow 1 + \mathfrak{m}$ and $b : \underline{I} \times \{\varsigma\} \longrightarrow 1 + \mathfrak{m}$ be multipliable in $1 + \mathfrak{m}$. Write $d := a \sqcup \{(\varsigma, 1 + \varepsilon)\} \sqcup b : I_\varsigma \longrightarrow 1 + \mathfrak{m}$. Assume that for each $i \in \underline{I}$, we have*

$$\varepsilon \cdot a(i) = b(i) \cdot a(i) \cdot \varepsilon. \tag{5.6}$$

*Then $\varepsilon \cdot (\Pi_I\, a) = \mathrm{ev}_d(M_I) \cdot (\Pi_I\, a) \cdot \varepsilon$.*

**Proof.** We have a strongly linear evaluation morphism $\mathrm{ev}_d : k\langle\langle I_\varsigma \rangle\rangle \longrightarrow A$. Let $\mathfrak{p} := \mathrm{Ker}(\mathrm{ev}_d)$, so $\mathfrak{p}$ is a summability ideal of $k\langle\langle I_\varsigma \rangle\rangle$. In view of (5.6), we have $\mathcal{R}(\underline{I}) \subseteq \mathfrak{p}$, so $\mathfrak{q}(\underline{I}) \subseteq \mathfrak{p}$. Let $P := M_I \cdot \mathrm{OP}_I \cdot X_\varsigma - X_\varsigma \cdot \mathrm{OP}_I$. We have $P \in \mathfrak{q}(\underline{I})$ by Lemma 5.1, so $\mathrm{ev}_d(P) = 0$, whence

$$\mathrm{ev}_d(M_I) \cdot (\Pi_I\, a) \cdot \varepsilon = \mathrm{ev}_d(M_I \cdot \mathrm{OP}_I \cdot X_\varsigma) = \mathrm{ev}_d(X_\varsigma \cdot \mathrm{OP}_I) = \varepsilon \cdot (\Pi_I\, a).$$ □

## 5.2 Application to Noetherian series

In the sequel, we assume that $k$ has characteristic zero. Let $(M, \cdot, 1, <)$ be an ordered monoid. Let $\mathbb{A} = k((M))$ be a summability algebra of Noetherian series over $k$ as per [4, Section 3]. So $\mathbb{A}$ is the $k$-algebra, under the convolution product, of functions $a : M \longrightarrow k$ whose support $\operatorname{supp} a = \{m \in M : a(m) \neq 0\}$ is a Noetherian subset of $(M, <)$, i.e. a well-founded subset without infinite antichains.

For $a, b \in \mathbb{A}$, we write $a \prec b$ if $b \neq 0$ and for all $m \in \operatorname{supp} a$, there is an $m' \in \operatorname{supp} b$ with $m' < m$. A $k$-linear function $\phi : \mathbb{A} \longrightarrow \mathbb{A}$ is said *contracting* if for all $a \in \mathbb{A} \setminus \{0\}$, we have $\phi(a) \prec a$. We write $\mathrm{Lin}^+_\prec(\mathbb{A})$ for the set of strongly linear and contracting functions $\mathbb{A} \longrightarrow \mathbb{A}$. The summability algebra $k\,\mathrm{Id}_{\mathbb{A}} + \mathrm{Lin}^+_\prec(\mathbb{A}) \subseteq \mathrm{Lin}^+(\mathbb{A})$ under composition has evaluations [4, Theorem 3.17]. It follows by Theorem 4.15 that $(\mathrm{Id}_{\mathbb{A}} + \mathrm{Lin}^+_\prec(\mathbb{A}), \circ)$ has a structure of multipliability group. Furthermore, writing $\mathbb{A}^\prec := \{a \in \mathbb{A} : a \prec 1\}$, the subalgebra $k + \mathbb{A}^\prec$ of $\mathbb{A}$ is a summability subalgebra that has evaluations.

**Theorem 5.3.** *The subset $1\text{-}\mathrm{Aut}^+_k(\mathbb{A}) := (\mathrm{Id}_{\mathbb{A}} + \mathrm{Lin}^+_\prec(\mathbb{A})) \cap \mathrm{End}^+(\mathbb{A})$ is a multipliability subgroup of $(\mathrm{Id}_{\mathbb{A}} + \mathrm{Lin}^+_\prec(\mathbb{A}), \circ)$.*

**Proof.** By [4, Theorem 3.17], the set $\log(1\text{-}\mathrm{Aut}^+_k(\mathbb{A}))$ is that of strongly linear and contracting derivations on $\mathbb{A}$. This is a summability Lie subalgebra of $\mathrm{Lin}^+_\prec(\mathbb{A})$ by [4, Proposition 2.2], so Proposition 4.16 applies and yields the result. □

**Corollary 5.4.** *The group $1\text{-}\mathrm{Aut}^+_k(\mathbb{A})$ has a natural structure of multipliability group.*

***Proof.*** *Since $1\text{-}\mathrm{Aut}^+_k(\mathbb{A})$ is a multipliability subgroup of $\mathrm{Id}_{\mathbb{A}} + \mathrm{Lin}^+_\prec(\mathbb{A})$, the structure of multipliability group on $\mathrm{Id}_{\mathbb{A}} + \mathrm{Lin}^+_\prec(\mathbb{A})$ induces a structure of full multipliability group on $1\text{-}\mathrm{Aut}^+_k(\mathbb{A})$.* □

We denote this multipliability structure by $\Pi$. Given $a \in \mathbb{A}$, we write $\mu_a$ for the product with $a$ on the left, i.e.

$$\begin{array}{rcl}\mu_a : \mathbb{A} & \longrightarrow & \mathbb{A} \\ b & \longmapsto & a \cdot b.\end{array}$$

If $\varepsilon \prec 1$, then $\mu_{1+\varepsilon}$ lies in $\mathrm{Id}_{\mathbb{A}} + \mathrm{Lin}^+_{\prec}(\mathbb{A})$. We write $\mu_{1+\mathbb{A}^{\prec}}$ for the subgroup of $\mathrm{Id}_{\mathbb{A}} + \mathrm{Lin}^+_{\prec}(\mathbb{A})$ of functions $\mu_{1+\varepsilon}$ where $\varepsilon$ ranges in $\mathbb{A}^{\prec}$.

**Lemma 5.5.** *The set $\mu_{1+\mathbb{A}^{\prec}}$ is a multipliability subgroup of $\mathrm{Id}_{\mathbb{A}} + \mathrm{Lin}^+_{\prec}(\mathbb{A})$ and*

$$\begin{array}{rcl}1+\mathbb{A}^{\prec} & \longrightarrow & \mu_{1+\mathbb{A}^{\prec}} \\ 1+\varepsilon & \longmapsto & \mu_{1+\varepsilon}\end{array}$$

*is a strongly multiplicative isomorphism.*

**Proof.** Let $J = (\underline{J}, <) \in \mathbf{Los}$ and let $(\mu_{1+\varepsilon_j})_{j \in \underline{J}}$ be $J$-multipliable in $\mathrm{Id}_{\mathbb{A}} + \mathrm{Lin}^+_{\prec}(\mathbb{A})$. So $(\mu_{1+\varepsilon_j} - \mathrm{Id}_{\mathbb{A}})_{j \in \underline{J}} = (\mu_{\varepsilon_j})_{j \in \underline{J}}$ is summable in $\mathrm{Lin}^+(\mathbb{A})$. This means in particular that the family $((\mu_{1+\varepsilon_j} - \mathrm{Id}_{\mathbb{A}})(1))_{j \in \underline{J}} = (\varepsilon_j)_{j \in \underline{J}}$ is summable in $\mathbb{A}$, whence in turn $(1+\varepsilon_j)_{j \in \underline{J}}$ is $J$-multipliable in $1+\mathbb{A}^{\prec}$. Write $1+\varepsilon$ for its product. For $a \in \mathbb{A}$, we have

$$\begin{aligned}\Big(\prod_J (\mu_{1+\varepsilon_j})_{j \in \underline{J}}\Big)(a) &= (\mathrm{ev}_{\mu_{1+\varepsilon_j} - \mathrm{Id}_{\mathbb{A}}}(\mathrm{OP}_J))(a) \\ &= \left(\sum_{\theta = (\theta_1, \ldots, \theta_n) \in J^\star} \mu_{\varepsilon_{\theta_1}} \circ \cdots \circ \mu_{\varepsilon_{\theta_n}}\right)(a) \\ &= \sum_{\theta = (\theta_1, \ldots, \theta_n) \in J^\star} \varepsilon_{\theta_1} \cdots \varepsilon_{\theta_n} \cdot a \\ &= \left(\sum_{\theta = (\theta_1, \ldots, \theta_n) \in J^\star} \varepsilon_{\theta_1} \cdots \varepsilon_{\theta_n}\right) \cdot a \\ &= (\mathrm{ev}_{(\varepsilon_j)_{j \in \underline{J}}}(\mathrm{OP}_J)) \cdot a \\ &= (1+\varepsilon) \cdot a.\end{aligned}$$

This shows that $\prod_J (\mu_{1+\varepsilon_j})_{j \in \underline{J}} = \mu_{1+\varepsilon}$, so $\mu_{1+\mathbb{A}^{\prec}}$ is a multipliability subgroup and $1+\delta \mapsto \mu_{1+\delta}$ is strongly multiplicative. □

If $\phi : \mathbb{A} \longrightarrow \mathbb{A}$ is a function, then we write $1\text{-}\mathrm{Aut}^+_{\phi}(\mathbb{A})$ for the subset of $1\text{-}\mathrm{Aut}^+_k(\mathbb{A})$ of automorphisms $\sigma$ that satisfy a chain rule

$$\phi \circ \sigma = \mu \circ \sigma \circ \phi$$

for a $\mu \in \mu_{1+\mathbb{A}^{\prec}}$. Since $\sigma \circ \mu_{1+\varepsilon} \circ \sigma^{-1} = \mu_{\sigma(1+\varepsilon)} \in \mu_{1+\mathbb{A}^{\prec}}$ for all $\varepsilon \in \mathbb{A}^{\prec}$, and $\sigma \in 1\text{-}\mathrm{Aut}^+_k(\mathbb{A})$, the set $1\text{-}\mathrm{Aut}^+_{\phi}(\mathbb{A})$ is a subgroup of $1\text{-}\mathrm{Aut}^+_k(\mathbb{A})$.

**Theorem 5.6.** *Assume that $\phi \in \mathrm{Lin}^+_{\prec}(\mathbb{A})$ and assume that there is a $\xi \in \mathbb{A}$ with $\phi(\xi) = 1$ and $\phi(a) \prec 1$ whenever $a \prec \xi$. Then $1\text{-}\mathrm{Aut}^+_{\phi}(\mathbb{A})$ is a multipliability subgroup of $1\text{-}\mathrm{Aut}^+_k(\mathbb{A})$.*

**Proof.** We fix an $\xi \in \mathbb{A}$ as in the statement. Let $I = (\underline{I}, <)$ be a linearly ordered set and let $(\sigma_i)_{i \in \underline{I}}$ be family in $1\text{-}\mathrm{Aut}^+_{\phi}(\mathbb{A})$ which is multipliable in $1\text{-}\mathrm{Aut}^+_k(\mathbb{A})$ with product $\sigma$, and let $\mu^{\sigma_i} \in \mu_{1+\mathbb{A}^{\prec}}$ with $\mu^{\sigma_i} \circ \sigma_i \circ \phi = \phi \circ \sigma$ for each $i \in \underline{I}$. The family $(\sigma_i - \mathrm{Id}_{\mathbb{A}})_{i \in \underline{I}}$ is summable in $\mathrm{Lin}^+(\mathbb{A})$, so $(\delta_i)_{i \in \underline{I}} := (\sigma_i(\xi) - \xi)_{i \in \underline{I}}$ is summable in $\mathbb{A}$. By strong linearity of $\phi$, so is the family $(\varepsilon_i)_{i \in \underline{I}} := (\phi(\delta_i))_{i \in \underline{I}}$. Let $i \in \underline{I}$. We have $\delta_i \prec \xi$, whence $\varepsilon_i \in \mathbb{A}^{\prec}$. Furthermore

$$1 + \varepsilon_i = 1 + \phi(\delta_i) = \phi(\sigma_i(\xi)) = \mu^{\sigma_i}(\sigma(\phi(\xi))) = \mu^{\sigma_i}(\sigma(1)) = \mu^{\sigma_i}(1).$$

We deduce since $(1+\varepsilon_i)_{i\in\underline{I}}$ is multipliable that $(\mu^{\sigma_i})_{i\in\underline{I}}$ is multipliable.

Consider the families $a:\underline{I}\longrightarrow 1\text{-Aut}^+(\mathbb{A})$ and $b:\underline{I}\times\{\varsigma\}\longrightarrow 1\text{-Aut}^+(\mathbb{A})$ given by

$$a(i):=\sigma_i \qquad \text{and} \qquad b(i,\varsigma):=\mu^{\sigma_i}$$

for each $i\in\underline{I}$. Also write $d:=a\sqcup b\sqcup\{(\varsigma,\mathrm{Id}_{\mathbb{A}}+\phi)\}$. Note that $(a,b)$ satisfies (5.6) with respect to $\phi$, so Lemma 5.2 gives

$$\phi\circ\sigma=\mathrm{ev}_d(M_I)\circ\sigma\circ\phi.$$

It suffices to show in order to conclude that $\mathrm{ev}_d(M_I)\in\mu_{1+\mathbb{A}^{\prec}}$. Recall that

$$\mathrm{ev}_d(M_I)=\mathrm{ev}_d(\mathrm{ev}_{Q_I-1}(\mathrm{OP}_I))=\mathrm{ev}_{(\mathrm{ev}_d(Q_I(i)-1))_{i\in\underline{I}}}(\mathrm{OP}_I)=\prod_I\,(\mathrm{ev}_d(Q_I(i)))_{i\in\underline{I}}$$

is an ordered product by Proposition 3.4. In view of Lemma 5.5, it suffices to show that each $\mathrm{ev}_d(Q_I(i))$ for $i\in\underline{I}$ lies in $\mu_{1+\mathbb{A}^{\prec}}$. Let $i\in\underline{I}$. We have $Q_I(i)=\mathrm{OP}_{i^-}\cdot(1+X_{(i,\varsigma)})\cdot\mathrm{OP}_{i^-}^{-1}$, where $\mathrm{ev}_d(1+X_{(i,\varsigma)})=\mu^{\sigma_i}\in\mu_{1+\mathbb{A}^{\prec}}$. Now $\varphi:=\mathrm{ev}_d(\mathrm{OP}_{i^-})$ lies in $1\text{-Aut}_k^+(\mathbb{A})$, so $\mathrm{ev}_d(Q_I(i))=\varphi\circ\mu^{\sigma_i}\circ\varphi^{\mathrm{inv}}$ lies in $\mu_{1+\mathbb{A}^{\prec}}$. □

# 6 Applications to flat transseries

We now turn to the field $\mathbb{T}$ of log-exp transseries of [5, 7, 6]. This is a direct limit of real algebras of Noetherian series, whose elements $s$ can be seen as formal series $s:\mathfrak{T}\longrightarrow\mathbb{R}$ with reverse well-ordered support, where $(\mathfrak{T},\cdot,1,\prec)$ is a totally ordered Abelian group called the group of transmonomials. We represent an $s\in\mathbb{T}$ as a formal sum $s=\sum_{\mathfrak{m}\in\mathfrak{M}}s(\mathfrak{m})\,\mathfrak{m}$, thus identifying $\mathfrak{T}$ with a subgroup of $\mathbb{T}^{\times}$. The function

$$\begin{aligned}\mathfrak{d}:\mathbb{T}^{\times}&\longrightarrow\mathfrak{T}\\ s&\longmapsto\min\operatorname{supp}s\end{aligned}$$

is a valuation on $\mathbb{T}$ with value group $(\mathfrak{T},\cdot,1,\succ)$ (note the reversed ordering). For $s,t\in\mathbb{T}^{\times}$, we write $s\asymp t$ if $\mathfrak{d}(s)=\mathfrak{d}(t)$, $s\preccurlyeq t$ if $\mathfrak{d}(s)\preccurlyeq\mathfrak{d}(t)$ and $s\prec t$ if $\mathfrak{d}(s)\prec\mathfrak{d}(t)$. We also have a lexicographic ordering $<$ on $\mathbb{T}$ for which it is an ordered field: a series $s\in\mathbb{S}^{\times}$ is positive if $s(\mathfrak{d}(s))>0$ as a real number. We write $\mathbb{T}^{>}=\{s\in\mathbb{S}:s>0\}$. There is an isomorphism $\log:(\mathbb{T}^{>},\cdot,1,<)\longrightarrow(\mathbb{S},+,0,<)$ whose inverse is denoted exp.

We write $x$ for the identity series in $\mathbb{T}$, which is a transmonomial. Roughly speaking, elements of $\mathbb{T}$ can be obtained from $x$ by combining exponentials, logarithms, products and infinite sums with real coefficients.

Consider the sugroup $\mathfrak{M}\subseteq\mathfrak{T}$ of elements $\mathfrak{m}$ such that there exists an $n\in\mathbb{N}$ with $x^{-n}\prec\mathfrak{m}\prec x^n$. For instance $\log x$, $\mathrm{e}^{\sqrt{\log x}}$ and $x^r$ are in $\mathfrak{M}$ for all $r\in\mathbb{R}$, but $\mathrm{e}^{-x},\mathrm{e}^{(\log x)^2}\notin\mathfrak{M}$. We define $\mathbb{S}$ to be the subfield $\mathbb{S}=\{s\in\mathbb{T}:\operatorname{supp}s\subseteq\mathfrak{M}\}$, called the field of *flat transseries* (see [1, Section 3]). It contains all logarithmic transseries. We write

$$\begin{aligned}\mathbb{S}^{>}&=\{s\in\mathbb{S}:s>0\},\\ \mathbb{S}^{>\mathbb{R}}&:=\{s\in\mathbb{S}:\forall s\in\mathbb{R},s>r\}=\{s\in\mathbb{S}:s>0\wedge s\succ 1\}\quad\text{and}\\ \mathcal{T}&:=\{s\in\mathbb{S}:s-x\prec x\}=x+x\,\mathbb{S}^{\prec}.\end{aligned}$$

## 6.1 Composition and derivation

We recall important facts on $\mathbb{S}$ and its structure (see [1, Section 3.1]).

i. The field $\mathbb{S}$ has a natural summability structure [13, 4, 8], and for each $s\in\mathbb{S}$, the family $(s(\mathfrak{m})\,\mathfrak{m})_{\mathfrak{m}\in\mathfrak{M}}$ is summable with sum $s$.

ii. There is a unique strongly linear derivation $\partial : \mathbb{S} \longrightarrow \mathbb{S}$ on $\mathbb{S}$ with $\partial(x) = 1$ and $\partial(\log(s)) = \frac{\partial(s)}{s}$ for all $s \in \mathbb{S}^{>}$.

iii. For all $s \in \mathbb{S}^{\times}$, we have $\partial(s) \preccurlyeq x^{-1} s$. In particular, for all $t \prec x$, the function $\mu_t \circ \partial$ is contracting.

iv. For each $t \in \mathbb{S}^{>\mathbb{R}}$, there is a unique strongly linear and strictly increasing morphism of algebras $\triangle_t : \mathbb{S} \longrightarrow \mathbb{S}$ with

$$\triangle_t(x) = t \quad \text{and} \tag{6.1}$$

$$\forall s \in \mathbb{S}^{>}, \triangle_t(\log s) = \log \triangle(s). \tag{6.2}$$

v. Write $s \circ t := \triangle_t(s)$ for all $s \in \mathbb{S}$, $t \in \mathbb{S}^{>\mathbb{R}}$. For $s \in \mathbb{S}$ and $t, u \in \mathbb{S}^{>\mathbb{R}}$, we have $t \circ u \in \mathbb{S}^{>\mathbb{R}}$ and

$$s \circ (t \circ u) = (s \circ t) \circ u.$$

vi. For each $s \in \mathbb{S}$ and $t \in \mathcal{T}$, the family $(\frac{\partial^{[i]}(s)}{i!} (t - x)^i)_{k \in \mathbb{N}}$ is summable with $s \succ \partial(s) (t - x) \succ \cdots \succ \frac{\partial^{[i]}(s)}{i!} (t - x)^i$ and

$$\triangle_t(s) = \sum_{i \in \mathbb{N}} \frac{\partial^{[i]}(s)}{i!} (t - x)^i.$$

In particular $\triangle_t \in 1\text{-Aut}_{\mathbb{R}}^{+}(\mathbb{S})$.

For $s, t \in \mathbb{S}^{>\mathbb{R}}$, we have $\triangle_t(\triangle_s(x)) = s \circ t$ whence $\triangle_t \circ \triangle_s = \triangle_{s \circ t}$ by iv. In particular $\mathcal{T}$ is a group under $\circ$ and we have a group isomorphism

$$\begin{array}{rcl} \mathcal{T} & \longrightarrow & \triangle_{\mathcal{T}} := \{\triangle_t : t \in \mathcal{T}\} \\ t & \longmapsto & \triangle_t^{\text{inv}}. \end{array}$$

## 6.2 Residues

We define the *residue* $\operatorname{res}(f)$ of an $s \in \mathcal{T} \setminus \{x\}$ to be the unique term $r\,\mathfrak{m} \in \mathbb{R}^{\times} \mathfrak{M}^{\prec x}$ with $s - x - r\,\mathfrak{m} \prec \mathfrak{m}$. So $s - x - \operatorname{res}(s) \prec \operatorname{res}(s)$. We also set $\operatorname{res}(x) := 0$.

For $s, t \in \mathcal{T} \setminus \{x\}$, writing $s = x + r_s \mathfrak{m}_s + \delta_s$ and $t = x + r_t \mathfrak{m}_t + \delta_t$ where $\delta_s \prec \operatorname{res}(s) = r_s \mathfrak{m}_s$ and $\delta_t \prec \operatorname{res}(t) = r_t \mathfrak{m}_t$, we have

$$\begin{aligned} s \circ t &= t + r_s \mathfrak{m}_s \circ t + \delta_s \circ t \\ &= t + r_s \mathfrak{m}_s \circ (x + r_t \mathfrak{m}_t + \delta_t) + \delta_s \circ t \\ &= x + r_t \mathfrak{m}_t + r_s \mathfrak{m}_s + (r_s \mathfrak{m}_s \circ (x + r_t \mathfrak{m}_t + \delta_t) - r_s \mathfrak{m}_s) + \delta_t + \delta_s \circ t. \end{aligned}$$

We have $\delta_t \prec \operatorname{res}(t)$ by definition, and $\delta_s \circ t \sim \delta_s \prec \operatorname{res}(s)$ and $r_s \mathfrak{m}_s \circ (x + r_t \mathfrak{m}_t + \delta_t) - r_s \mathfrak{m}_s \prec r_s \mathfrak{m}_s$ by vi. From this, we deduce that:

$$(\operatorname{res}(s \circ t) = \operatorname{res}(s) \text{ if } \operatorname{res}(s) \succ \operatorname{res}(t)) \text{ and } (\operatorname{res}(s \circ t) = \operatorname{res}(t) \text{ if } \operatorname{res}(s) \prec \operatorname{res}(t)) \tag{6.3}$$

$$(\operatorname{res}(s \circ t) = \operatorname{res}(s) + \operatorname{res}(t)) \quad \text{if} \quad (\operatorname{res}(s) \asymp \operatorname{res}(t) \text{ and } \operatorname{res}(s) + \operatorname{res}(t) \neq 0) \tag{6.4}$$

$$(\operatorname{res}(s \circ t) \prec \operatorname{res}(s), \operatorname{res}(t)) \qquad \text{if} \qquad \operatorname{res}(s) + \operatorname{res}(t) = 0. \tag{6.5}$$

**Remark 6.1.** The structure $(\mathcal{T}, \circ, x, <)$ is a growth order group as per [3, Definition 2.27], and our notion of residue corresponds to that of [3, Definition 2.18].

## 6.3 The group of right compositions

In view of [1, Theorem 3.1], the set $\triangle_{\mathcal{T}}$ coincides with $\{\exp(f\partial) : f \in \mathbb{S} \wedge f \prec x\}$. From this observation, we deduce the following

**Lemma 6.2.** *Let $s\in\mathcal{T}$ and $r\in\mathbb{R}$. Then $\triangle_s^{[r]}=\triangle_{s^{[r]}}$ for a unique series $s^{[r]}\in\mathcal{T}$.*

**Proof.** Write $\triangle_s=\exp(f\partial)$ for an $f\prec x$. Then $\triangle_s^{[r]}=\exp(r f\partial)\in\triangle_{\mathcal{T}}$, hence the result. □

We also deduce that for $s\in\mathcal{T}$, wee have $s^{[r+r']}=s^{[r]}\circ s^{[r']}$ for all $r,r'\in\mathbb{R}$, and $\{s^{[r]}:r\in\mathbb{R}\}$ is the centralizer of $s$ in $\mathcal{T}$ if $s\neq x$. Let $s\in\mathcal{T}\setminus\{x\}$ and $r\in\mathbb{R}^{\times}$. We have $s^{[r]}=\triangle_s^{[r]}(x)=x+r\,\phi(x)+\binom{r}{2}\phi^{[2]}(x)+\cdots$ where $\phi$ is contracting and $\phi(x)=\sum_{i>0}\frac{\partial^{[i]}(x)}{i!}(s-x)^i=s-x$. So $(s^{[r]}-x)-r\,(s-x)\prec x-s\asymp\operatorname{res}(s)$, i.e. $(s^{[r]}-x)-r\operatorname{res}(s)\prec\operatorname{res}(s)$, whence

$$\operatorname{res}(s^{[r]})=r\operatorname{res}(s). \tag{6.6}$$

**Proposition 6.3.** *We have $\triangle_{\mathcal{T}}=1\text{-}\mathrm{Aut}_{\partial}^{+}(\mathbb{S})$.*

**Proof.** For $f\in\mathbb{S}$ with $f\prec x$, the automorphism $\exp(f\partial)$ is a power series evauated at $f\partial$, so it commutes with $f\partial$. In particular, it satisfies a chain rule with respect to $\partial$. So $\triangle_{\mathcal{T}}\subseteq 1\text{-}\mathrm{Aut}_{\partial}^{+}(\mathbb{S})$.

Conversely $\sigma\in 1\text{-}\mathrm{Aut}_{\mathbb{R}}^{+}(\mathbb{S})$ and write $\sigma=\mathrm{Id}_{\mathbb{S}}+\phi$ where $\phi:\mathbb{S}\longrightarrow\mathbb{S}$ is contracting. Suppose that $\sigma$ satisfies a chain rule $\sigma(s)'=f\sigma(s')$ for all $s\in\mathbb{S}$. We claim that $\sigma$ commutes with log on $\mathbb{S}^{>}$. Note by (6.2) that it suffices to show the result for the automorphism $\triangle_{x+x^{-1}}\circ\sigma$. Since $\triangle_{x+x^{-1}}$ satisfies a chain rule, we may replace $\sigma$ with $\triangle_{x+x^{-1}}\circ\sigma$. In view of vi, the map $\triangle_{x+x^{-1}}\circ\sigma-\mathrm{Id}_{\mathbb{S}}$ sends $\mathbb{S}$ to $\mathbb{S}^{\prec}$, so we may assume that $\phi(s)\prec 1$ for all $s\in\mathbb{S}$. Let $s\in\mathbb{S}^{>}$. We have

$$(\sigma(\log s))'=f\sigma((\log s)')=f\sigma\left(\frac{s'}{s}\right)=\frac{f\sigma(s')}{\sigma(s)}=\frac{\sigma(s)'}{\sigma(s)}=(\log\sigma(s))'.$$

So $(\sigma(\log s))'=(\log\sigma(s))'$, whence there is an $r\in\mathbb{R}$ with $\sigma(\log s)=\log\sigma(s)+r$. We have

$$\sigma(\log s)=\log s+\phi(\log s)=\log(s+\phi(s))+r=\log s+\varepsilon(s)+r$$

where $\varepsilon:=\sum_{i>0}\frac{\partial^{[i]}(\log)(s)}{i!}\phi(s)^i$ and $\phi(\log s)$ are both infinitesimal. It follows that $r=\phi(\log s)-\varepsilon$ is infinitesimal, so $r=0$. Thus $\sigma$ comutes with log on $\mathbb{S}^{>}$. We deduce that $\sigma\in\triangle_{\mathbb{S}^{>\mathbb{R}}}$ whence $\sigma\in\triangle_{\mathcal{T}}$ as $\sigma(x)\sim x$ by definition of $1\text{-}\mathrm{Aut}_{\mathbb{R}}^{+}(\mathbb{S})$. □

**Theorem 6.4.** *The set $\triangle_{\mathcal{T}}=\{\triangle_f:f\in\mathcal{T}\}$ is a multipliability subgroup of $1\text{-}\mathrm{Aut}_{\mathbb{R}}^{+}(\mathbb{S})$.*

**Proof.** The set $\triangle_{\mathcal{T}}$ coincides with $1\text{-}\mathrm{Aut}_{\partial}^{+}(\mathbb{S})$. Now $\partial$ satisfies the hypotheses of Theorem 5.6, so $1\text{-}\mathrm{Aut}_{\partial}^{+}(\mathbb{S})$ is a multipliability subgroup of $1\text{-}\mathrm{Aut}_{\mathbb{R}}^{+}(\mathbb{S})$. □

We can also directly obtain the result by applying Proposition 4.16 to $\exp(\mathbb{S}^{\prec x}\partial)$, however the formalism of chain rules is more general, as

**Corollary 6.5.** *The group $\triangle_{\mathcal{T}}$ has a natural structure of multipliability group.*

**Proof.** This follows from Theorem 6.4 and Proposition 2.7. □

## 6.4 Decompositions along scales

In view of Theorem 6.4, we have a natural and full multipliability structure $\Pi$ on $\triangle_{\mathcal{T}}$, which the structure on $1\text{-}\mathrm{Aut}_{\mathbb{R}}^{+}(\mathbb{S})$ extends. Recall that we have an isomorphism

$$\begin{aligned}\mathcal{E}:(\triangle_{\mathcal{T}},\circ)&\longrightarrow(\mathcal{T},\circ)\\ \sigma&\longmapsto\sigma^{\mathrm{inv}}(x).\end{aligned}$$

By Corollary 2.10, we obtain a full multipliability structure $\Pi^* := \Pi^{\mathcal{E}}$ on $\mathcal{T}$ given by

$$\begin{aligned}
\operatorname{dom} \Pi_I^* &= \{\mathcal{E} \circ \sigma : \sigma \in \operatorname{dom} \Pi_I\} \\
&= \{(s_i)_{i \in \underline{I}} : (\triangle_{s_i}^{\mathrm{inv}})_{i \in \underline{I}} \in \operatorname{dom} \Pi_I\} \\
&= \{(s_i)_{i \in \underline{I}} : (\triangle_{s_i})_{i \in \underline{I}} \in \operatorname{dom} \Pi_{I^*}\} && \text{(by MG7)} \\
& \quad \{(s_i)_{i \in \underline{I}} : (\triangle_{s_i})_{i \in \underline{I}} \in \operatorname{dom} \Pi_I\}, \quad \text{and} && \text{(by Fullness)} \\
\Pi_I^* (s_i)_{i \in \underline{I}} &:= \mathcal{E}(\Pi_I ((\triangle_{s_i}^{\mathrm{inv}})_{i \in \underline{I}})) \\
&= (\Pi_I ((\triangle_{s_i}^{\mathrm{inv}})_{i \in \underline{I}}))^{\mathrm{inv}}(x) \\
&= (\Pi_{I^*} (\triangle_{s_i})_{i \in \underline{I}})(x); && \text{(by MG7)}
\end{aligned}$$

for all $(s_i)_{i \in \underline{I}} \in \operatorname{dom} \Pi_I^*$.

**Proposition 6.6.** *Let $I$ be a set and let $(s_i)_{i \in I} \in \mathcal{T}^I$ be a family. Then $(s_i)_{i \in I}$ is multipliable in $(\mathcal{T}, \Pi^*)$ if and only if $(s_i - x)_{i \in I}$ is summable in $(\mathbb{S}, \Sigma)$.*

**Proof.** Write $\sigma_i := \triangle_{s_i}$, $\phi_i := \sigma_i - \mathrm{Id}_{\mathbb{S}}$ and $\delta_i := s_i - x \prec x$ for each $i \in I$. Recall that by definition, the family $(s_i)_{i \in I}$ is multipliable in $(\mathcal{T}, \Pi^*)$ if and only if $(\phi_i)_{i \in I}$ is Lin-summable in $\mathrm{Lin}^+(\mathbb{S})$. For $i \in I$, and $s \in \mathbb{S}$, we have

$$\phi_p(s) = \sum_{p > 0} \frac{s^{(p)}}{p!} \delta_i^p = \left( \sum_{i > 0} \frac{1}{p!} \mu_{\delta_i}^{[p]} \circ \partial^{[p]} \right)(s).$$

Assume that $(s_i)_{i \in I}$ is multipliable. Then $(\phi_i(x))_{i \in I} = (\delta_i)_{i \in I}$ is summable in $\mathbb{S}$.

Conversely, assume that $(\delta_i)_{i \in I}$ is summable in $\mathbb{S}$. We claim that the family $(\phi_i)_{i \in I}$ is Lin-summable in $\mathrm{Lin}_{\prec}^+(\mathbb{S})$. Indeed, for $i \in I$, we have $\delta_i \prec x$ and $\delta_i' \prec 1$, so the functions $\mu_{\delta_i} \circ \partial$ and $\mu_{\delta_i'}$ are contracting. Furthermore, since $\partial$ is strongly linear, the family $(\delta_i')_{i \in I}$ is summable. So both $(\mu_{\delta_i})_{i \in I}$ and $(\mu_{\delta_i'})_{i \in I}$ are Lin-summable, whence also $(\mu_{\delta_i} \circ \partial)_{i \in I}$ is Lin-summable.

Set $J := (I \times \{0, 1\}) \sqcup \{0\}$, and consider the family $\phi : J \longrightarrow \mathrm{Lin}_{\prec}^+(\mathbb{S})$ with $\phi(0) = \partial$, and $\phi(i, 0) = \mu_{\delta_i} \circ \partial$ and $\phi(i, 1) = \mu_{\delta_i'}$ for each $i \in I$. By SS4 in $\mathrm{Lin}^+(\mathbb{S})$, the family $\phi$ is Lin-summable. Note that $\partial \circ \mu_{\delta_i} = \phi(i, 0) + \phi(i, 1)$ for all $i \in I$. An easy induction on $p > 0$, shows that there are an $n_p > 0$, integers $m_{p,1}, \ldots, m_{p,n_p}$ and words $\theta_{p,1}(i), \ldots, \theta_{p,n_p}(i) \in ((\{i\} \times \{0, 1\}) \sqcup \{0\})^{\star}$ of length $p$ such that

$$\mu_{\delta_i}^{[p]} \circ \partial^{[p]} = \mathrm{ev}_{\phi} \left( \sum_{\ell=1}^{n_p} m_{p,\ell} X_{\theta_{p,\ell}(i)} \right). \tag{6.7}$$

Note that since each $\theta_{p,\ell}(i)$ for $\ell \in \{1, \ldots, n_p\}$ has length $p$, the sum

$$P := \sum_{i \in I} \sum_{p > 0} \sum_{\ell=1}^{n_p} \frac{1}{p!} m_{p,\ell} X_{\theta_{p,\ell}(i)} \in \mathbb{R}\langle\!\langle J \rangle\!\rangle$$

is well-defined. The algebra $\mathbb{R}\,\mathrm{Id}_{\mathbb{S}} + \mathrm{Lin}_{\prec}^+(\mathbb{S})$ has evaluations, so the family $(P(\theta)\ \mathrm{ev}_{\phi}(X_\theta))_{\theta = (\theta_1, \ldots, \theta_n) \in J^{\star}}$ is Lin-summable. By SS3 in $\mathrm{Lin}^+(\mathbb{S})$, so is $(\sum_{\ell=1}^{n_p} P(\theta_{p,\ell}(i))\, \mathrm{ev}_{\phi}(X_{\theta_{p,\ell}(i)}))_{p > 0 \wedge i \in I} = (\frac{1}{p!} \mu_{\delta_i}^{[p]} \circ \partial^{[p]})_{p > 0 \wedge i \in I}$. We deduce again with SS3 that $(\sum_{p > 0} \frac{1}{p!} \mu_{\delta_i}^{[p]} \circ \partial^{[p]})_{i \in I} = (\phi_i)_{i \in I}$ is Lin-summable. Therefore $(\sigma_i)_{i \in I}$ is multipliable in $1\text{-}\mathrm{Aut}_{\mathbb{R}}^+(\mathbb{S})$, so $(s_i)_{i \in I}$ is multipliable in $\mathcal{T}$. □

We write $\mathfrak{M}^{\prec x} = \{\mathfrak{m} \in \mathfrak{M} : \mathfrak{m} \prec x\}$.

**Definition 6.7.** *A* **scale** *on $\mathcal{T}$ is a function $S : \mathfrak{M}^{\prec x} \times \mathbb{R} \longrightarrow \mathcal{T}$ such that*

*i. $S(\mathfrak{m}, 0) = x$ for all $\mathfrak{m} \in \mathfrak{M}^{\prec x}$.*

*ii.* $\mathrm{res}(S(\mathfrak{m},r)) = r\,\mathfrak{m}$ *for all* $(\mathfrak{m},r) \in \mathfrak{M}^{\prec x} \times \mathbb{R}$.

*iii. for all reverse well-ordered subsets* $\mathfrak{W} \subseteq \mathfrak{M}^{\prec x}$ *and all families* $(r_w)_{\mathfrak{w}\in\mathfrak{W}} \in \mathbb{R}^{\mathfrak{W}}$*, the family* $(S(\mathfrak{w},r_{\mathfrak{w}}) - x)_{\mathfrak{w}\in\mathfrak{W}}$ *is summable in* $\mathbb{S}$.

**Example 6.8.** The functions $S_0(\mathfrak{m},r) := x + r\,\mathfrak{m}$ and $S_1(\mathfrak{m},r) := (x+\mathfrak{m})^{[r]}$ are scales on $\mathbb{S}$. For $S_0$, this is clear. Let us turn to $S_1$. For $(\mathfrak{m},r) \in \mathfrak{M}^{\prec x} \times \mathbb{R}$, we have $\mathrm{res}(S_1(\mathfrak{m},r)) = r\,\mathfrak{m}$ by (6.6) and $S_1(\mathfrak{m},0) = x$ by definition. Let $\mathfrak{W} \subseteq \mathfrak{M}^{\prec x}$ be reverse well-ordered and let $(r_{\mathfrak{w}})_{\mathfrak{w}\in\mathfrak{W}} \in \mathbb{R}^{\mathfrak{W}}$. Then $(S(\mathfrak{w},1) - x)_{\mathfrak{w}\in\mathfrak{W}} = (\mathfrak{w})_{\mathfrak{w}\in\mathfrak{W}}$ is summable, so by Proposition 6.6, the family $(\triangle_{S(\mathfrak{w},1)})_{\mathfrak{w}\in\mathfrak{W}}$ is multipliable. Therefore $(\phi_{\mathfrak{w}})_{\mathfrak{w}\in\mathfrak{W}} := (\triangle_{S_1(\mathfrak{w},1)} - \mathrm{Id}_{\mathbb{S}})_{\mathfrak{w}\in\mathfrak{W}}$ is Lin-summable in $\mathrm{Lin}^+_{\prec}(\mathbb{S})$. Since $\mathbb{R}\,\mathrm{Id}_{\mathbb{S}} + \mathrm{Lin}^+_{\prec}(\mathbb{S})$ has evaluations, so is $(\binom{r_{\mathfrak{w}}}{j}\,\phi_{\mathfrak{w}}^{[i]})_{\mathfrak{w}\in\mathfrak{W}\wedge i\in\mathbb{N}^{>0}}$. By **SS3** in $\mathbb{R}\,\mathrm{Id}_{\mathbb{S}} + \mathrm{Lin}^+_{\prec}(\mathbb{S})$, the family $(\sum_{i>0}\binom{r_{\mathfrak{w}}}{i}\,\phi_{\mathfrak{w}}^{[i]})_{\mathfrak{w}\in\mathfrak{W}} = (\triangle_{S(\mathfrak{w},1)}^{[r_{\mathfrak{w}}]} - \mathrm{Id}_{\mathbb{S}})_{\mathfrak{w}\in\mathfrak{W}}$ is also Lin-summable. Therefore $(\triangle_{S_1(\mathfrak{w},1)}^{[r_{\mathfrak{w}}]}(x) - x)_{\mathfrak{w}\in\mathfrak{W}} = (S_1(\mathfrak{w},r_{\mathfrak{w}}))_{\mathfrak{w}\in\mathfrak{W}}$ is summable in $\mathbb{S}$.

As a corollary of Proposition 6.6, we obtain:

**Corollary 6.9.** *Let $S$ be a scale on $\mathbb{S}$. For $r : \mathfrak{M}^{\prec x} \longrightarrow \mathbb{R}$ with reverse well-ordered support, the family $(S(\mathfrak{w},r(\mathfrak{w})))_{\mathfrak{w}\prec x}$ is multipliable in $(\mathcal{T},\Pi^*)$.*

**Lemma 6.10.** *Let $S$ be a scale on $\mathbb{S}$ and let $\lhd$ be a linear ordering on $\mathfrak{M}^{\prec x}$. For $r : \mathfrak{M}^{\prec x} \longrightarrow \mathbb{R}$ with reverse well-ordered support in $(\mathfrak{M}^{\prec x}, <)$, writing $\mathfrak{w}_0 := \min \operatorname{supp} r \in \mathfrak{M}^{\prec x}$ we have $\mathrm{res}(\Pi^*_{(\mathfrak{M}^{\prec x},\lhd)}(S(\mathfrak{w},r(\mathfrak{w})))_{\mathfrak{w}\prec x}) = r(\mathfrak{w}_0)\,\mathfrak{w}_0$.*

**Proof.** Write $I := (\mathfrak{M}^{\prec x}, \lhd)$ and $s := \Pi^*_{(\mathfrak{M}^{\prec x},\prec)}\,(S(\mathfrak{w}, r(\mathfrak{w})))_{\mathfrak{w}\prec x}$. We first prove that $\mathrm{res}(s) \preccurlyeq \mathfrak{w}_0$. As in Corollary 6.9, writing $\sigma_{\mathfrak{w}} := \triangle_{S(\mathfrak{w},r(\mathfrak{w}))}$ and $\phi_{\mathfrak{w}} := \sigma_{\mathfrak{w}} - \mathrm{Id}_{\mathbb{S}}$ for each $\mathfrak{w} \in \mathfrak{M}^{\prec x}$, we have $\phi_{\mathfrak{w}} = \sum_{p>0} \frac{1}{p!}\,\mu^{[p]}_{S(\mathfrak{w},r(\mathfrak{w}))-x} \circ \partial^{[p]}$ and $s = (\mathrm{ev}_{(\phi_{\mathfrak{w}})_{\mathfrak{w}\prec x}}(\mathrm{OP}_{I^*}))(x)$. Since each $\phi_{\mathfrak{w}}$ is contracting, we have $\mathrm{ev}_{(\phi_{\mathfrak{w}})_{\mathfrak{w}\prec x}}(X_\theta)(x) \prec \mathrm{ev}_{(\phi_{\mathfrak{w}})_{\mathfrak{w}\prec x}}(X_{\theta_n})(x)$ for all $\theta = (\theta_1,\ldots,\theta_n) \in (I^*)^\star$ with $n > 1$. Since $\{(j) : j \in \underline{I}\} \subseteq (I^*)^\star$, and $\phi_{\mathfrak{w}}(x) \prec \phi_{\mathfrak{y}}(x)$ whenever $\mathfrak{w} \prec \mathfrak{y}$ in $\mathfrak{M}^{\prec x}$, it follows that $\mathrm{res}(s) = \mathrm{res}(x + \phi_{\mathfrak{w}_0}(x))$. Now $\phi_{\mathfrak{w}_0}(x) = (S(\mathfrak{w}_0, r(\mathfrak{w}_0)) - x) \cdot 1 + 0 = S(\mathfrak{w}_0, r(\mathfrak{w}_0)) - x$, so $\mathrm{res}(s) = \mathrm{res}(S(\mathfrak{w}_0, c(\mathfrak{w}_0))) = r(\mathfrak{w}_0)\,\mathfrak{w}_0$. □

Let $\kappa$ be a cardinal number that does not embed into $(\mathfrak{M}^{\prec x}, \succ)$. Let $N : \kappa \longrightarrow \{-1, 1\}$ be a function and consider the linear ordering $<_N$ on $\kappa$ of Section 1.3. If $\lambda \in \kappa$, $\mathfrak{w} : (\lambda, \in) \longrightarrow (-\infty, 1)$ is a strictly decreasing function and $r : \lambda \longrightarrow \mathbb{R}$ is a function, then by Corollary 6.9 and **MG3**, the family $s := (S(\mathfrak{w}(\gamma), r(\gamma)))_{\gamma\in\lambda}$ is multipliable in $(\mathcal{T},\Pi^*)$. We define

$$N_S(\mathfrak{w},r) := \Pi^*_{(\lambda,<_N)}\, s \in \mathcal{T}.$$

Writing

$$\begin{aligned} N_{S,\gamma}(\mathfrak{w},r) &:= \Pi^*_{([\gamma,\lambda),<_N)}\,(s \upharpoonright [\gamma,\lambda)), \\ N^L_{S,\gamma}(\mathfrak{w},r) &:= \Pi^*_{(L_N(\gamma,\lambda),<_N)}\,(s \upharpoonright L_N(\gamma,\lambda)) \quad \text{and} \\ N^R_{S,\gamma}(\mathfrak{w},r) &:= \Pi^*_{(R_N(\gamma,\lambda),<_N)}\,(s \upharpoonright R_N(\gamma,\lambda)), \end{aligned}$$

the axiom **MG4c** gives

$$N_S(\mathfrak{w},r) = N^L_{S,\gamma}(\mathfrak{w},r) \circ N_{S,\gamma}(\mathfrak{w},r) \circ N^R_{S,\gamma}(\mathfrak{w},r). \tag{6.8}$$

In view of Lemma 1.2 and **MG4c**, given $\eta < \gamma$, we can further expand

$$N^L_{S,\gamma}(\mathfrak{w},r) = N^L_{S,\eta}(\mathfrak{w},r) \circ t^L_{\gamma,\eta}, \qquad \text{and} \qquad N^R_{S,\gamma}(\mathfrak{w},r) = t^R_{\gamma,\eta} \circ N^R_{S,\eta}(\mathfrak{w},r), \tag{6.9}$$

where

$$t^L_{\gamma,\eta} := \Pi^*_{(L_n(\gamma,\lambda)\setminus L_N(\eta,\lambda),<_N)}(s\upharpoonright L_N(\gamma,\lambda)\setminus L_N(\eta,\lambda)) \quad \text{and}$$
$$t^R_{\gamma,\eta} := \Pi^*_{(R_N(\gamma,\lambda)\setminus R_N(\eta,\lambda),<_N)}(s\upharpoonright R_N(\gamma,\lambda)\setminus R_N(\eta,\lambda)).$$

By Lemma 6.10, we have

$$\operatorname{res}(t^L_{\gamma,\eta}), \operatorname{res}(t^R_{\gamma,\eta}) \prec \mathfrak{w}(\eta). \tag{6.10}$$

**Theorem 6.11.** *Let $S$ be a scale on $\mathbb{S}$ and let $N:\kappa\longrightarrow\{-1,1\}$ be a function. For each $s\in\mathcal{T}$, there are a unique ordinal $\lambda_0<\kappa$, a unique strictly decreasing function $\mathfrak{w}:\lambda_0\longrightarrow(-\infty,1)$ and a unique function $r:\lambda_0\longrightarrow\mathbb{R}^\times$ such that $s=N_S(\mathfrak{w},r)$.*

**Proof.** Let us first prove the existence. We may assume that $s\neq x$, otherwise the empty family yields the result. We define a strictly increasing function $\mathfrak{w}:\lambda\longrightarrow(-\infty,1)$ and a map $r:\lambda\longrightarrow\mathbb{R}^\times$ on an ordinal $\lambda<\kappa$ as follows. Let $\lambda<\kappa$ such that for all $\gamma<\lambda$, the maps $\mathfrak{w}$ and $r$ are defined on $\gamma$ and $\mathfrak{w}$ is strictly decreasing on $\gamma$, and we have

$$\operatorname{res}((N^L_{S,\eta}(\mathfrak{w}\upharpoonright\gamma, r\upharpoonright\gamma))^{[-1]}\circ s\circ(N^R_{S,\eta}(\mathfrak{w}\upharpoonright\gamma,r\upharpoonright\gamma))^{[-1]})\preccurlyeq\mathfrak{w}(\eta) \tag{6.11}$$

for all $\eta<\gamma$. If $\lambda$ is a limit, then $\mathfrak{w}$ and $r$ are defined on $\lambda$ and $\mathfrak{w}$ is strictly decreasing on $\lambda$. The condition (6.11) is trivially satisfied for $\lambda$.

Assume that $\lambda=\gamma+1$ is a successor. Let $\mathfrak{w}_0, r_0$ be the extensions of $\mathfrak{w}, r$ to $\lambda$ with $(\mathfrak{w}_0(\gamma), r_0(\gamma))=(\mathfrak{w}(\gamma)-1,0)$. Set

$$s_1 := (N^L_{S,\gamma}(\mathfrak{w}_0,r_0))^{[-1]}\circ s\circ(N^R_{S,\gamma}(\mathfrak{w}_0,r_0))^{[-1]}.$$

If $s_1=x$, then we set $\lambda=\gamma$, and we have

$$s=N^L_{S,\gamma}(\mathfrak{w}_0,r_0)\circ N^R_{S,\gamma}(\mathfrak{w}_0,r_0)=N_{S,\gamma}(\mathfrak{w},r). \tag{6.12}$$

Assume that $s_1\neq x$. For $\eta<\gamma$, we can write

$$N^L_{S,\gamma}(\mathfrak{w}_0,r_0)=N^L_{S,\eta}(\mathfrak{w}_0,r_0)\circ t^L_{\gamma,\eta} \qquad \text{and} \qquad N^R_{S,\eta}(\mathfrak{w}_0,r_0)=t^R_{\gamma,\eta}\circ N^R_{S,\gamma}(\mathfrak{w}_0,r_0)$$

as in (6.9), where $\operatorname{res}(t^L_{\gamma,\eta}), \operatorname{res}(t^R_{\gamma,\eta})\prec x^{e(\eta)}$. Note that $s_1=(t^L_{\gamma,\eta})^{[-1]}\circ t\circ(t^R_{\gamma,\eta})^{[-1]}$, where

$$t:=(N^L_{S,\eta}(\mathfrak{w}_0,r_0))^{[-1]}\circ s\circ(N^R_{S,\eta}(\mathfrak{w}_0,r_0))^{[-1]}$$

has residue $\operatorname{res}(t)\preccurlyeq\mathfrak{w}(\eta)$ by the induction hypothesis. We deduce with (6.3–6.5) that $\operatorname{res}(s_1)\preccurlyeq\mathfrak{w}(\eta)$. Set $r\,\mathfrak{w}:=\operatorname{res}(s_1)$. Assume that $N(\gamma)=1$. and write $r(\gamma)\,\mathfrak{w}(\gamma):=\operatorname{res}(s_1\circ S(r,\mathfrak{w})^{[-1]})$. We have $\mathfrak{w}(\gamma)\prec\mathfrak{w}(\eta)$ for all $\eta<\gamma$ by (6.5), so this extends the functions $\mathfrak{w}, r$ as desired. Now we have $(N^L_{S,\gamma}(\mathfrak{w},r))^{[-1]}\circ s\circ(N^R_{S,\gamma}(\mathfrak{w},r))^{[-1]}=s_1\circ S(\mathfrak{w},r)^{[-1]}$, so the sequences $\mathfrak{w}, r$ satisfy (6.11) for $\lambda$. If $N(\gamma)=-1$, then define $r(\gamma)\,\mathfrak{w}(\gamma)=\operatorname{res}(S(\mathfrak{w},r)^{[-1]}\circ s_1)$. Symmetric arguments show that this extends the functions as required.

Since $(\mathfrak{w}(\gamma))_{\gamma<\lambda}$ is a strictly decreasing sequence in $(\mathfrak{M}^{\prec x},\prec)$, the inductive process must stop at an ordinal $\lambda_0<\kappa$, and we have $s=N_S(\mathfrak{w},r)$ as in (6.12) for that ordinal.

We next prove the unicity. Let $\lambda_1<\kappa$, let $\mathfrak{w}_1:\lambda_1\longrightarrow\mathfrak{M}^{\prec x}$ be a strictly decreasing function and let $r_1:\lambda_1\longrightarrow\mathbb{R}^\times$ be a function. Without loss of generality, we may assume that $\lambda\leqslant\lambda_1$. Assume that $(\mathfrak{w},r)\neq(\mathfrak{w}_1,r_1)$ and let $\alpha\in\lambda\cup\lambda_1$ be minimal to witness this. That is, either $\alpha=\lambda<\lambda_1$, or $\alpha<\lambda$ and $(\mathfrak{w}(\alpha),r(\alpha))\neq(\mathfrak{w}_1(\alpha),r_1(\alpha))$. Writing

$$N_S(\mathfrak{w}_1,r_1)=N^L_{S,\alpha}(\mathfrak{w}_1,r_1)\circ N_{S,\alpha}(\mathfrak{w}_1,r_1)\circ N^R_{S,\alpha}(\mathfrak{w}_1,r_1),$$

we see that $N_S(\mathfrak{w}, r) = N^L_{S,\alpha}(\mathfrak{w}_1, r_1) \circ d \circ N^R_{S,\alpha}(\mathfrak{w}_1, r_1)$ for $d = \prod_{([\alpha,\lambda),<_N)} (S(\mathfrak{w}(\gamma), r(\gamma)))_{\gamma\in[\alpha,\lambda)}$. We have $d = x$ if $\alpha = \lambda$, and $\operatorname{res}(d) = r(\alpha)\, \mathfrak{w}(\alpha)$ by Lemma 6.10 if $\alpha < \lambda$. In the first case, since $\operatorname{res}(N_{S,\alpha}(\mathfrak{w}_1, r_1)) = r_1(\alpha)\, \mathfrak{w}_1(\alpha) \neq 0$, we have $N_{S,\alpha}(\mathfrak{w}_1, r_1) \neq x$ so $N_S(\mathfrak{w}_1, r_1) \neq N_S(\mathfrak{w}, r)$. In the second case, since since $\operatorname{res}(N_{S,\alpha}(\mathfrak{w}_1, r_1)) = r_1(\alpha)\, \mathfrak{w}_1(\alpha) \neq r(\alpha)\, \mathfrak{w}(\alpha)$, we have $N_S(\mathfrak{w}_1, r_1) \neq N_S(\mathfrak{w}, r)$. □

**Theorem 6.12.** *For each $a \in \mathcal{T}$, there are a unique ordinal $\lambda$, a unique strictly decreasing function $s: \lambda \longrightarrow \{x + \mathfrak{m} : \mathfrak{m} \prec x\}$ and a unique function $r: \lambda \longrightarrow \mathbb{R}^\times$ with*

$$a = \cdots \circ s_\gamma^{[r(\gamma)]} \circ \cdots \circ s_1^{[r(1)]} \circ s_0^{[r(0)]}.$$

**Proof.** Take $S$ to be the scale $S_1$ of Example 6.8, and take $N$ to be the constant function $N(\gamma) = -1$ for all $\gamma < \omega_1$. Then apply Theorem 6.11. □

Taking $S$ to be the scale $S_1$ of Example 6.8, Theorem 6.11 gives a representation of $\mathcal{T}$ as an infinite semidirect product $\mathcal{T} \simeq \rtimes_{\mathfrak{w}\in\mathfrak{M}^{\prec x}} (\mathbb{R}, +, 0)$ of copies of $(\mathbb{R}, +, 0)$, with the trivial multipliability structure, along the linear ordering $\mathfrak{M}^{\prec x}$. This gives an instance of a positive answer of [3, Question 1].

**Remark 6.13.** Let us give a brief explanation of how to extend this to the group $\mathbb{T}^{>\mathbb{R}}$ of positive infinite log-exp transseries. Any transseries $f \in \mathbb{T}^{>\mathbb{R}}$ can be written uniquely as $f = f_0 \circ \exp_l(x)$ where $l \in \mathbb{Z}$ is the "level" (see [20, Definition 1]) or "exponentiality" (see [7, (4.1.22)]) of $f$, and $\exp_l(x)$ is the $l$-fold iterate of exp evaluated at $x$. The transseries $f_0$ has level 0, which means that for large enough $n \in \mathbb{N}$, we have $\exp_{-n}(x) \circ f_0 \circ \exp_n(x) = x + \delta$ for a $\delta \prec x$. For $i > 0$, consider the subgroup $\mathfrak{T}_i = \{\mathfrak{m} \in \mathfrak{T} : \log \mathfrak{m} \prec \exp_{i-1}(x)\}$ of $\mathfrak{T}$, write $\mathbb{T}_i := \mathbb{R}[\![\mathfrak{T}_i]\!]$ and set

$$\mathcal{T}_i := \{x + \varepsilon \in x + \mathbb{T}_i : \varepsilon \prec 1 \wedge \log|\varepsilon| \succcurlyeq \exp_{i-1}(x)\}.$$

We claim that $x + \delta$ may be written uniquely as a finite composition

$$x + \delta = (x + \delta_d) \circ \cdots \circ (x + \delta_0),$$

where $x + \delta_0 \in \mathcal{T}$ and $x + \delta_i \in \mathcal{T}_i$ for $i \in \{1, \ldots, d\}$. More precisely, the group $\mathcal{T}_{\mathbb{T}} = \{x + u \in \mathbb{T} : u \prec x\}$ is an iterated semidirect product of the $\mathcal{T}_d$'s. Each $\mathcal{T}_i$ has the form

$$\mathcal{T}_i = \exp\left(\mathbb{T}_i^{\prec}\left(\prod_{j<i} \exp_j(x)\right)\partial\right)(x)$$

and can be represented as $1\text{-Aut}^+_\partial(\mathbb{T}_i)$, thus elements in $\mathcal{T}_i$ can be written uniquely as infinite compositions of elements of the form $(x + \mathfrak{m})^{[r]}$, $\mathfrak{m} \in \mathfrak{T}_i, r \in \mathbb{R}^\times$. Using a natural construction of semidirect products of mulgipliability groups, one recovers an expression of $f$ by expanding $(x + \delta)$ as an infinite composition, and expanding $f = \exp_n(x) \circ (x + \delta) \circ \exp_{-n}(x) \circ \exp_l(x)$ by "distributing" the conjugacy map $g \mapsto \exp_n(x) \circ g \circ \exp_{-n}(x)$ over th infinite product, as prescribed by **MG6**.